\documentclass{article}

\usepackage{amsmath, amsthm, amscd, amsfonts, amssymb, graphicx, color}
\usepackage[bookmarksnumbered, colorlinks, plainpages]{hyperref}
\usepackage{pspicture,pst-all,multido}
\theoremstyle{definition}

\theoremstyle{remark}

\numberwithin{equation}{section}

\newcommand{\SMALL}{\scriptsize}

\newcommand{\C}{\mathbb{C}}
\newcommand{\R}{\mathbb{R}}
\newcommand{\N}{\mathbb{N}}

\renewcommand{\ker}{\operatorname{ker}}
\newcommand{\F}{\mathcal{F}}

\newcommand{\dom}{\operatorname{dom}}
\newcommand{\ran}{\operatorname{ran}}

\newcommand{\norm}[1]{\parallel\!\!#1\!\!\parallel}
\newcommand{\sgn}{\operatorname{sgn}}
\renewcommand{\i}{\operatorname{i}}
\renewcommand{\d}{\operatorname{d}}
\renewcommand{\deg}{\operatorname{deg}}
\newcommand{\e}{\operatorname{e}}

\newcounter{envcount}%

\newenvironment{Def}%
{\vspace{\bigskipamount}\refstepcounter{envcount}\textbf{(\theenvcount)\enspace Definition.}}%
  {\vspace{\bigskipamount}}

{\vspace{\bigskipamount}\refstepcounter{envcount}\textbf{(\theenvcount)\enspace Summary.}}%
  {\vspace{\bigskipamount}}

\newenvironment{HS}%
{\vspace{\bigskipamount}\refstepcounter{envcount}\textbf{(\theenvcount)\enspace Hardy spaces.}}%
  {\vspace{\bigskipamount}}

{\vspace{\bigskipamount}\refstepcounter{envcount}\textbf{(\theenvcount)\enspace Modulus of a Hardy function.}}%
  {\vspace{\bigskipamount}}

\newenvironment{RTO}%
{\vspace{\bigskipamount}\refstepcounter{envcount}\textbf{(\theenvcount)\enspace Remark on Toeplitz Operators.}}%
  {\vspace{\bigskipamount}}

\newenvironment{SBS}%
{\vspace{\bigskipamount}\refstepcounter{envcount}\textbf{(\theenvcount)\enspace Semibounded symbol.}}%
  {\vspace{\bigskipamount}}

{\vspace{\bigskipamount}\refstepcounter{envcount}\textbf{(\theenvcount)\enspace Bounded nonnegative symbol.}}%
  {\vspace{\bigskipamount}}

\newenvironment{Exa}%
{\vspace{\bigskipamount}\refstepcounter{envcount}\textbf{(\theenvcount)\enspace Example.}}%
  {\vspace{\bigskipamount}}

{\vspace{\bigskipamount}\refstepcounter{envcount}\textbf{(\theenvcount)\enspace Example}}%
  {\vspace{\bigskipamount}}

{\vspace{\bigskipamount}\enspace \textbf{Definition.}}%
  {\vspace{\bigskipamount}}

\newenvironment{The}%
{\vspace{\bigskipamount}\refstepcounter{envcount}\textbf{(\theenvcount)\enspace Theorem.}\itshape}%
  {\vspace{\bigskipamount}\upshape}

{\vspace{\bigskipamount}\refstepcounter{envcount}\textbf{(\theenvcount)\enspace Theorem}\itshape}%
  {\vspace{\bigskipamount}\upshape}

{\vspace{\bigskipamount}\refstepcounter{envcount}\textbf{(\theenvcount)\enspace Proposition.}\itshape}%
  {\vspace{\bigskipamount}\upshape}

\newenvironment{Cor}%
{\vspace{\bigskipamount}\refstepcounter{envcount}\textbf{(\theenvcount)\enspace Corollary.}\itshape}%
  {\vspace{\bigskipamount}\upshape}

\newenvironment{Lem}%
{\vspace{\bigskipamount}\refstepcounter{envcount}\textbf{(\theenvcount)\enspace Lemma.}\itshape}%
  {\vspace{\bigskipamount}\upshape}

{\vspace{\bigskipamount}\refstepcounter{envcount}\textbf{(\theenvcount)\enspace}}%
  {\vspace{\bigskipamount}}

\theoremstyle{definition}
\swapnumbers

\begin{document}
\setcounter{page}{1}
\pagenumbering{arabic}

\begin{center}
{\Large  Unbounded Wiener-Hopf Operators and Isomorphic Singular Integral Operators}\\  

\vspace{1cm}
Domenico P.L. Castrigiano\\
Technischen Universit\"at M\"unchen, Fakult\"at f\"ur Mathematik, M\"unchen, Germany\\

\smallskip

{\it E-mail address}: {\tt
castrig\,\textrm{@}\,ma.tum.de
}
\end{center}

\begin{quote} Some basics of a theory of unbounded Wiener-Hopf operators (WH), not existing so far, are developed. 
The alternative is shown that the domain of a WH is either zero or dense. The symbols with  non-trivial WH, which are called proper, are determined explicitly by an integrability property. WH are characterized by shift invariance.  We study in detail WH with rational symbols showing that they are  densely defined  closed and have finite dimensional kernels and deficiency spaces. The latter spaces as well as the domains, ranges, spectral and Fredholm points are explicitly determined. Another topic concerns semibounded  WH.  There is a canonical representation of a semibounded WH using a product of a closable operator and its adjoint. The Friedrichs extension is obtained replacing the operator by its closure.  The polar decomposition  gives rise to a Hilbert space isomorphism relating a semibounded WH  to  a singular integral operator of  type  Hilbert transformation. This remarkable relationship, which allows to transfer results and methods reciprocally, is new also in the thoroughly studied case of bounded WH. \\


\end{quote}


\section{Introduction}

There is an increasing interest in unbounded Toeplitz and Toeplitz-like operators, which  will concern also the related  Wiener-Hopf operators (WH) (see (\ref{RTO})). So far, results on WH $W_\kappa$ with unbounded symbol $\kappa$ are scarce and probably in the literature  there exists  no introduction to this subject. So sec.\,\ref{UBWHO} deals with preliminaries and basics regarding unbounded WH. In particular we are  concerned with conditions on the symbol $\kappa$ ensuring that the domain of $W_\kappa$ is either the whole space or dense or   trivial, and prove  that $\dom W_\kappa$ is either trivial or dense. The symbols with  non-trivial WH, which are called proper, are determined by a useful integrability property.  A  classical result on the eigenvalues of a WH is shown to remain valid  in the unbounded case. It  implies that non-trivial symmetric WH have no eigenvalues. A further result characterizes WH by their invariance under unilateral shifts.  In sec.\,\ref{RS} WH with rational symbols are studied. They constitute a welcome source of densely defined closed operators with finite index.  An explicit description of the    domains, ranges, kernels, deficiency spaces, spectral and Fredholm points is given.  The remainder of this article deals in sec.\,\ref{SBWHO} and sec.\,\ref{ISIO} with densely defined semibounded WH. A semibounded operator $W_\kappa$ can be expressed by a  product of a closable operator  $A$ and its adjoint. Replacing $A$ by its closure one obtains quite naturally a self-adjoint extension  $\tilde{W}_\kappa$. It is proven to coincide   with the Friedrichs extension. Inverting the order of the factors one obtains a singular integral operator $L_\phi$ of type Hilbert transformation.  For the operators of the mention type being not trivial  there is a necessary condition analogous to that for WH. The  self-adjoint extensions  $\tilde{L}_\phi$ and $\tilde{W}_\kappa$ are isometric, which  follows from the polar decomposition of $\overline{A}$. Actually $\tilde{W}_\kappa$ is  Hilbert space isomorphic to the reduction of  $\tilde{L}_\phi$ on $\ker(\tilde{L}_\phi)^\perp$, and the spectral representations of $\tilde{L}_\phi$ and $\tilde{W}_\kappa$ can be achieved in an explicit manner from each other. It is worth noting that this relationship is new also for bounded WH thus contributing to the well-developed theory of the latter.
To conclude, this method is illustrated  by a non-trivial example  diagonalizing Lalescu's operator and the isometrically related singular integral operator. In  \cite[sec.\,3.3]{C19} the spectral representations of   $W_{1_{[-1,1]}}$ and  the finite Hilbert transformation are related to each other by this method.\\

\textbf{Notations.}  Let $\F$ denote the Fourier transformation on $L^2(\R)$. For   measurable $E\subset \R$ introduce the projection  $P_E:L^2(\R)\to L^2(E)$, $(P_Ef)(x):=f(x)$. (For convenience define $L^2(E)=\{0\}$, $P_E=0$ if $E=\emptyset$.) Its adjoint  $P_E^*$ is  the injection $(P_E^*f)(x)=f(x)$ for $x\in E$ and $=0$ otherwise. Note $$P_EP_E^*=I_{L^2(E)}, \quad P_E^*P_E=M(1_E)$$ 
with $M(1_E)$ the multiplication by the indicator function $1_E$ for $E$. We call $E$ \textbf{proper} if neither $E$ nor the complement $\R\setminus E$ is a null set. Put $\R_+:= ]0,\infty[$ and $P_+:=P_{\R_+}$. Analogously define $P_-$. \\
\hspace*{6mm}
Throughout  let $\kappa:\R\to\C$ denote a measurable function and $M(\kappa)$  the multiplication operator  by $\kappa$  in $L^2(\R)$ with dense domain $\{f\in L^2(\R): \kappa f\in L^2(\R)\}$.  $M(\kappa)$ is normal  satisfying    $M(\kappa)^*=
M(\overline{\kappa})$ and $M(\overline{\kappa})M(\kappa)=M(|\kappa|^2)$.  Moreover,  $M(\kappa)$ is self-adjoint if and only if $\kappa$ is almost real, and  $M(\kappa)$ is nonnegative if and only if $\kappa\ge 0$ a.e.

\begin{Def}\label{GWHO}  The operator in $L^2(\R_+)$
$$W_\kappa:=P_+\mathcal{F}M(\kappa)\mathcal{F}^{-1}P_+^*$$  is called the Wiener-Hopf operator  (WH) with symbol $\kappa$. Occasionally we write $W(\kappa)$ instead of $W_\kappa$. Clearly, $W_\kappa=W_{\kappa'}$ if $\kappa=\kappa'$ a.e.  Often we shall  refer to this tacitly. The symbol $\kappa$ is called \textbf{proper} if $\dom W_\kappa\ne \{0\}$.
\end{Def}\\
The theory of WH  with bounded symbol is well developed. We content ourselves to refer here to the book \cite[Chapter 9]{BS90} and to mention the origins \cite{WHP31}. Obviously in case of a bounded symbol   the operators $W_\kappa$ are bounded  with $\dom W_\kappa=L^2(\R_+)$  and   adjoint $W_\kappa^* = W_{\overline{\kappa}}$. $W_\kappa$ is the convolution on the real half line  with kernel $k$, i.e., 
$$(W_\kappa g)(x)=\int_0^\infty k(x-y)g(x)\operatorname{d}y$$
 if $\kappa\in L^\infty\cap L^2 $ and $k:=(2\pi)^{-1/2}\mathcal{F}\kappa$, or if
 $\kappa=\int\operatorname{e}^{\operatorname{i}(\cdot) y}k(y)\operatorname{d}y$ for  $k\in L^1(\R)$. For the case of integrable kernel there is the rather exhaustive theory by M. G. Krein \cite{K62}.
 Generally  the tempered distribution $k:=(2\pi)^{-1/2}\mathcal{F}\kappa$, where  $\kappa\in L^\infty(\R)$ is considered as a regular tempered distribution,   satisfies $\mathcal{F}M(\kappa)\mathcal{F}^{-1}u=k\star u$ for every Schwartz function  $u$ in the distributional sense (e.g.\,\cite[Theorem IX.4]{RS75}). 
   For instance the kernel for $W_{-\operatorname{sgn}}$ is the tempered distribution $k(x)=-\frac{1}{x}$ or that for $W_{-\tanh}$ equals $k(x)=\big(2\operatorname{i}\sinh(\pi x/2)\big)^{-1}$. 
 In the literature the generalizations of WH  stay mostly within the realm of bounded operators. One deals with  the traces (compressions) of bounded bijective operators in  Banach space on a closed subspace  \cite{DSS69}. 
The results concern the solvability of the associated Wiener-Hopf equations.

\section{Unbounded Wiener-Hopf Operators}\label{UBWHO}

As put it by  \cite{Y17} results on unbounded WH are practically inexistent. Indeed  they are scarce. See \cite[1.3]{Y18} for some notes.  
An important result 
is due to M. Rosenblum  \cite{R62}, \cite{R65}, obtained for Toeplitz operators 
  and hence valid for the  Hilbert space isomorphic WH $W_\kappa$ (see (\ref{RTO})). So in the case that the symbol $\kappa$ is real bounded below not almost constant and $(1+x^2)^{-1}\kappa$ is integrable,   \cite{R65} furnishes  the spectral representation of the extension
 $\tilde{W}_\kappa$, which is shown to be the Friedrichs extension  (\ref{RSBWHO}),\,(\ref{SBS}) and which  by  \cite{R65} is  absolutely continuous. --- There are investigations on unbounded general WH dealing with conditions for their invertibility \cite{Pe71}. ---
In \cite{Y17} real bounded below  Wiener-Hopf quadratic forms from distributional kernels $k$ are considered, and it is shown that such a  form determines a WH if and only if the form is closable or, equivalently, if and only if $\sqrt{2\pi}k$ is the Fourier transform of a locally integrable bounded below function $\kappa$ with integrable $(1+x^2)^{-n}\kappa$ for some $n\in\N$. Clearly $\kappa$ is the symbol and $\dom W_\kappa\supset C_c^\infty(\R_+)$ holds. See further \cite{Y18}. --- Furthermore, the methods applied for the study of unbounded analytic Toeplitz operators \cite{S08} can also produce results on unbounded WH, as   (\ref{DBWHO})(i), (\ref{AL2}), (\ref{RSBWHO})(c).\\

 Starting the preliminary remarks and anticipating briefly some results on WH,   
 note  first that  $\dom W_\kappa=\dom M(\kappa)\F^{-1}P_+^*$. Therefore $\kappa$ is proper if and only if  $\kappa h$ is square-integrable for some Hardy function (\ref{HS}) $h\in  \ran(\F^{-1}P_+^*)$. Clearly $\dom W_\kappa\subset \dom W_{\kappa'}$ if $|\kappa'|\le |\kappa|$. Moreover, $\dom W_\kappa=\dom W_{\overline{\kappa}}=\dom W_{|\kappa|}=\dom W(1+|\kappa|)$. Note
$\langle g',W_\kappa\, g\rangle = \langle W_{\overline{\kappa}}\, g', g\rangle$  for $g, g'\in\dom W_\kappa$. Hence, if $W_\kappa$ is densely defined then $W_\kappa\subset W_{\overline{\kappa}}^*$, whence $W_\kappa$ is closable and 
$\dom W_\kappa\subset \dom W_\kappa^*$ holds.  But unbounded WH may and may not be closed (\ref{RSKCK}),\,(\ref{ENCWH}). \\
\hspace*{6mm}
If $W_\kappa$ is densely defined then $W_\kappa$ is symmetric, i.e.,  $W_\kappa\subset W_\kappa^*$, if and only if $\kappa$ is almost real (\ref{DBWHO})(n). If $W_\kappa$ is densely defined symmetric, then $W_\kappa$ is bounded below if and only if $\kappa$ is essentially bounded below (\ref{DBWHO})(o). Recall that the numerical range $\{\langle g,W_\kappa g\rangle: g\in\dom W_\kappa,\;||g||=1\}$ of $W_\kappa$ is convex. (Indeed, the numerical range of every operator in Hilbert space is convex \cite{G70}.) It is determined in (\ref{DBWHO})(o),(p).\\
\hspace*{6mm}
If the symbol $\kappa$ is unbounded, then $\dom W_\kappa\ne L^2(\R_+)$ (\ref{DBWHO})(a).  The alternative holds that either $\dom W_\kappa$ is trivial or $\dom W_\kappa$ is dense. In other words, as shown in  (\ref{DBWHO})(i), if $\kappa$ is proper, then $W_\kappa$ is densely defined. In (\ref{DBWHO})(e) and (\ref{DBWHO})(g) explicit characterizations of proper symbols are given. There is also the useful criterion in  (\ref{DBWHO})(b) for $\kappa$ to be proper.  So proper symbols may have polynomial growth and countably many singularities with integrable logarithm like as  $\exp|x|^\alpha$, $-1<\alpha<0$. It is  easy to give examples of  non-proper symbols (\ref{DBWHO})(h).

\begin{HS}\label{HS}
Recall the Hardy spaces $H_\pm:=\ran (\F^{-1}P_\pm^*)$. Obviously $L^2(\R)=H_+\,\textcircled{$\perp$}\, H_-$ and $h\in H_+$ $\Leftrightarrow$ $\overline{h}\in H_-$ as well $h\in H_+$ $\Leftrightarrow$ $\check{h}\in H_-$ with $\check{h}(x):=h(-x)$. We tacitly refer to the well-known Paley-Wiener Theorem characterizing the Fourier transforms of  $L^2$-functions vanishing on a half-axis, see e.g. \cite[Theorem 95]{T48}. In particular $h\in H_+$ if and only if there is a $\phi$ holomorphic on the upper  half-plane such that its partial maps $\phi_y(x):=\phi(x+\i y)$ for $y>0$ satisfy $\phi_y\in L^2(\R)$,   $\{\norm{\phi_y}: y>0\}$ bounded,  and $\phi_y\to h$ for $y\to 0$  in the mean and pointwise a.e. Actually $\phi$ converges to $h$ non-tangentially a.e., and $\norm{\phi_y}\uparrow\norm{h}$ for $y\downarrow 0$  (see e.g.\,\cite[III 3.3, II 2.6]{R05}).
Moreover, every $h\in H_\pm\setminus \{0\}$ vanishes only on a null set. Indeed, according to a Luzin-Privalov Theorem  \cite[IV 2.5]{Pr65} a meromorphic function on the upper or lower half-plane which takes non-tangential boundary values zero on a set of positive Lebesgue measure is zero. 
The former property is  also an immediate consequence of  the following result  on the modulus of a Hardy function:
\begin{equation}\label{MHF}
 \textit{For } f\in L^2(\R)\setminus\{0\} \textit{ there is } h\in H_+ \textit{ satisfying } |h|=|f| \textit{ if and only if } \,\textrm{\small{$\frac{\operatorname{ln} |f|}{1+x^2}$}}\in L^1(\R)
 \end{equation}
 One proves (\ref{MHF}) using  the outer function with prescribed modulus on the torus (\cite[Chap.\,3, Def.\,1.1 and   Prop.\,3.2]{A14}) and  the Hilbert space isomorphism $\Gamma$ in (\ref{RTO}). Recall  also \cite[Theorem XII ]{PW34} for $\ln|h|\,/(1+x^2)\in L^1$ if $h\in H_+$.
 \\
\hspace*{6mm}
Let $M_+(\kappa)$ denote  the trace of $M(\kappa)$ on $H_+$, i.e.
\begin{equation}\label{WHHS}
M_+(\kappa)=P_{H_+}M(\kappa)P^*_{H_+}
\end{equation}
 with $P_{H_+}f$ the orthogonal projection of $f\in L^2(\R)$ on $H_+$. Note that $P_{H_+}\F^{-1}P_+^*: L^2(\R_+)\to H_+$
 is a Hilbert space isomorphism with its inverse $P_+\F P_{H_+}^*$, by which $W_\kappa$ is Hilbert space  isomorphic to  $M_+(\kappa)$. Often it is convenient to deal with  $M_+(\kappa)$ in
place of $W_\kappa$. Note that $M_+(\kappa)$ extends
the multiplication operator in $H_+$
 $$M(\kappa, H_+)h:=\kappa\, h\; \textit{ with } \,\dom M(\kappa, H_+)=\{h\in H_+:\kappa h\in H_+\}$$
The latter is Hilbert space isomorphic to its counterpart  the so-called analytic Toeplitz operator with symbol $\omega:=\kappa\circ C^{-1}$, see (\ref{RTO}). Sometimes $M_+(\kappa)$ and $M(\kappa, H_+)$ coincide as for instance for rational symbols $\kappa$ holomorphic in the upper half-plane (\ref{AL1})(c). We do not put forward this topic.\\
 \hspace*{6mm}
Finally $H_+^\infty$ is the set of all measurable bounded $\alpha:\R\to \C$ such that there is a bounded holomorphic $A$ on the upper half-plane with the partial maps $A_y\to \alpha$ for $y\to 0$ pointwise a.e. Actually $A$ converges to $\alpha$ non-tangentially a.e. (see e.g.\,\cite[III 3.3, II 2.6]{R05}).
\end{HS}

\begin{RTO}\label{RTO} Let the torus $\mathbb{T}$ be endowed with the normalized Lebesgue measure. The Hardy space $H^2(\mathbb{T})$ is the subspace of $L^2(\mathbb{T})$ with orthonormal basis  
$e_n(w):=w^n$, $n\in\N_0$. 
Given a measurable $\omega: \mathbb{T}\to \C$ then, quite analogous to (\ref{WHHS}), the Toeplitz  operator $T_\omega$ with symbol $\omega$ is defined by $$T_\omega :=P_{H^2(\mathbb{T})}M(\omega)\,P^*_{H^2(\mathbb{T})}$$
Let $\Gamma:H^2(\mathbb{T}) \to H_+$, $(\Gamma u)(x):=\frac{\i}{\sqrt{\pi}(x+\i)}u\big(C(x)\big)$ be the Hilbert space isomorphism based on the Cayley transformation  $C(x):=\frac{x-\i}{x+\i}$. One has
\begin{equation}\label{THS}
T_\omega=\Gamma^{-1} M_+(\omega\circ C)\, \Gamma
\end{equation}
Obviously by this well-known relationship results and methods regarding Toeplitz operators may be transferred for the study of WH and vice versa.  For unbounded Toeplitz and Toeplitz-like operators see \cite{R62} and \cite{R65} for integrable symbols, \cite{S08} for analytic Toeplitz operators, \cite{Y18} for semibounded Toeplitz operators, \cite{GHJR18} for rational symbols.
\end{RTO}

\begin{The}\label{DBWHO} \emph{(a)} $\dom W_\kappa=L^2(\R_+)$ $\Rightarrow$ $\kappa$ is bounded  $\Rightarrow$ $\dom W_\kappa=L^2(\R_+)$ and  $W_\kappa$ is bounded.\\

\emph{(b)}  If $\kappa\, qs\in L^2(\R)$, where $q$  is a polynomial and $s$ is  the inverse Fourier transform of a Schwartz function with support in $[0,\infty[$, then $qs\in\dom W_\kappa$ and $\kappa$ is proper. --- 
Suppose that $\kappa\, qs\in L^2(\R)$ for a polynomial $q$ with only real zeros and for every Schwartz function $s$. Then $\dom W_\kappa\supset\{q(\operatorname{\i\frac{\d}{\d x}})\phi: \phi\in C_c^\infty(\R_+)\}$. \\  
  
\emph{(c)}  $W_\alpha (\dom W_\kappa)\subset \dom W_\kappa$ $\forall$  $\alpha\in H_+^\infty$ and, equivalently, $\alpha h \in \dom M_+(\kappa)$ $\forall$  $\alpha\in H_+^\infty$, $h\in \dom M_+(\kappa)$.\\

\emph{(d)} For all $\alpha\in H_+^\infty$ one has \,$\emph{(1)}$\;\;$W_\kappa W_\alpha=W_{\alpha \kappa}$ and 
 \,$\emph{(2)}$\;\;$W_\alpha^* W_\kappa\subset W_{\overline{\alpha} \kappa}$. ---
 For $b\ge 0$ let $S_b:L^2(\R)\to L^2(\R)$ be the unilateral translation $S_bg(x):=g(x-b)$ if $x>b$ and $:=0$ otherwise. Then one has  the translational  invariance $$W_\kappa=S_b^*W_\kappa S_b$$
  
\emph{(e)}  Let $p\in]0,\infty]$. Let $j$ denote a real-valued function such that $\frac{j}{1+x^2}$ is integrable.  Then
 $$\kappa \textrm{ proper } \Leftrightarrow\; \kappa \e^j\in L^p(\R) \textrm{ for some }  j$$
 
\emph{(f)} $\kappa$ is proper if and only if $\kappa^2$ is proper. More generally, let $r>0$ and let $\kappa_1$, $\kappa_2$ be two symbols satisfying $|\kappa_1|= |\kappa_2|^r$. Then $\kappa_1$ is proper if and only if 
$\kappa_2$ is proper. Finally, if $\kappa_1$ and $\kappa_2$ are proper symbols then so are $\kappa_1\kappa_2$ and $\kappa_1+\kappa_2$.\\

\emph{(g)} $\kappa$ is proper if and only if $\frac{\ln(1+|\kappa|)}{1+x^2}$ is integrable.\\

\emph{(h)} $\kappa$ is not proper if $|\kappa(x)|$ increases not less than exponentially  for $x\to\infty$ or $x\to -\infty$, i.e., if there are positive constants  $a$, $\delta$, $\lambda$ such that $|\kappa(x)|\ge \delta \e^{\lambda |x|}$ for $x\ge a$ or $x\le -a$. \\

\emph{(i)} If $\kappa$ is proper, then $\dom W_\kappa$ is dense. More precisely, if $h\in\dom M_+(\kappa)\setminus\{0\}$, then there is $h'\in H_+$ with $|h|=|h'|$ such that $H_+^\infty \frac{1}{x+\i}h'$ is dense and contained in $\dom M_+(\kappa)$.\\

\emph{(j)}  Put $\rho:= (1+|\kappa|)^{-1}$. Then $\kappa$ is proper if and only if $\overline{M(\rho)H_-}\cap M(\rho)H_+=\{0\}$ and $\kappa$ is not proper if and only if   $\overline{M(\rho)H_-}=L^2(\R)$.\\

\emph{(k)} Let $\kappa$  be not almost constant.  If  $\lambda\in\C$ is an eigenvalue of $W_\kappa$, then $\overline{\lambda}$ is not an eigenvalue of $W_{\overline{\kappa}}$. If  $\kappa$ is almost real, then $W_\kappa$ has no eigenvalues and in particular $W_\kappa$ is injective.\\

\emph{(l)}  Let $ \lambda\in\C$ such that $\kappa^{-1}(\{\lambda\})$ is proper.   Then $\lambda$ is not an eigenvalue of $W_\kappa$.\\

\emph{(m)}  Let $W_{\kappa_1}$ and $W_{\kappa_2}$ coincide on a dense set of $L^2(\R_+)$. Then $\kappa_1=\kappa_2$ a.e.\\

\emph{(n)} Let $W_\kappa$ be densely defined. Then $W_\kappa$ is symmetric if and only if $\kappa$ is almost real.\\

\emph{(o)} Let $W_\kappa$ be densely defined symmetric. Then $W_\kappa$ is  bounded below if and only if $\kappa$ is real  essentially bounded below, and the maximal lower bound  of  $W_\kappa$ equals the maximal essential lower bound of $\kappa$. If $W_\kappa$ is bounded below and not bounded, then the numerical range $\{\langle g,W_\kappa g\rangle: g\in\dom W_\kappa,\;||g||=1\}$ equals $]\alpha,\infty[$ with $\alpha$ the maximal lower bound.\\

\emph{(p)} Let $W_\kappa$ be bounded symmetric and not  a  multiple of $I$.  Then the numerical range 
equals $]a,b[$ with $a$ and $b$ the minimum and maximum, respectively, of the essential range of $\kappa$. 
\end{The}

{\it Proof.} (a) The second implication is obvious. As to the first one assume that $\kappa$ is not bounded. Then $f\in L^2(\R)$ exists with $\kappa f\notin L^2(\R)$. Write $f= h_++h_-$ with $h_\pm\in H_\pm$. As $\dom M_+(\kappa)=L^2(\R_+)$ one has $\kappa h_+\in L^2(\R)$, and hence  $\kappa h_-\notin L^2(\R)$. Similarly, $\kappa\overline{f}\notin L^2(\R)$, $\overline{f}=\overline{h}_-+\overline{h}_+$ with $\overline{h}_\mp\in H_\pm$, and 
$\kappa \overline{h}_-\in L^2(\R)$.  This contradicts $\kappa h_-\notin L^2(\R)$.\\
 (b) As to the first claim note that $s$ is a Schwartz function in $H_+$, whence $qs\in H_+$ by (\ref{AL1})(c).\\
\hspace*{6mm}
 Now let $\operatorname{D}$ denote the differential operator $\i\frac{\d}{\d x}$ and let $u$ be any Schwartz function with support in $\R_+$. Then $q(\operatorname{D})u$ is still such a function, and $\F^{-1}u$, $\F^{-1} \big(q(\operatorname{D})u\big)=q\F^{-1}u$ are Schwartz functions in $H_+$. Hence it remains to show that $\{q(\operatorname{D})\phi: \phi\in\mathcal{D}\}$ is dense in $L^2(\R_+)$, where $\mathcal{D}$ denotes the space of test functions   $C^\infty_c(\R_+)$. Assume $g\in L^2(\R_+)$ with $g \perp q(\operatorname{D})\phi$, i.e., $\int_{\R_+}\overline{g}\,q(\operatorname{D})\phi\d x=0$. The claim is $g=0$.\\
\hspace*{6mm} 
  Indeed, regarding $\overline{g}$ as a regular distribution in $\mathcal{D}'$, one has $\overline{g}\big(q(\operatorname{D})\phi\big)=q(-\operatorname{D})\overline{g}\,(\phi)=0$ for all test functions $\phi$, whence $q(-\operatorname{D})\overline{g}=0$. Thus $\overline{g}$ is a solution of the differential equation $q(-\operatorname{D})F=0$ for $F\in\mathcal{D}'$. As known all its solutions are regular.  Hence $\overline{g}\in T$, where $T$ denotes the space of linear combinations of functions on $\R_+$ of the kind $x\to x^k\e^{\i\lambda x}$, $k\in\N_0$ and $\lambda\in\R$. Then $G\in T$ for $G(x):=\int_0^x|g(t)|^2\d t$. One has $G(x)\to\, \norm{g}^2<\infty$ for $x\to \infty$. Write $G=\sum_{k=0}^n x^kA_k$ with $A_k$  a linear combination of periodic functions  $\e^{\i\lambda x}$. Then $A_n(x)\to 0$ for $x\to \infty$. Since $A_n$ is almost periodic this implies $A_n=0$. The result follows. \\
(c) Let $h\in\dom M_+(\kappa)$.  The claim  is  $M_+(\alpha)h\in\dom M_+(\kappa)$. Now $M_+(\alpha)h=\alpha h$ since obviously $\alpha h\in H_+$. Moreover $\kappa \alpha h\in L^2(\R)$, whence the claim.\\
(d) As to (1) note $h\in\dom \big(M_+(\kappa)M_+(\alpha)\big)$ $\Leftrightarrow$ $h\in H_+$, $P_{H_+}(\alpha h)\in \dom M_+(\kappa)$. Since $\alpha h\in H_+$ the latter is equivalent to $\alpha h\in\dom M_+(\kappa)$
$\Leftrightarrow$ $\kappa\alpha h\in L^2(\R)$ $\Leftrightarrow$ $h\in\dom  M_+(\alpha\kappa)$. --- Regarding (2) note $M^*_+(\alpha)=M_+(\overline{\alpha})$, $\dom \big(M^*_+(\alpha)  M_+(\kappa)\big)=\dom M_+(\kappa)$, and $P^*_{H_+} P_{H_+}=I-P^*_{H_-}P_{H_-}$ for $P_{H_-}:=I-P_{H_+}$. Hence for $h\in\dom M_+(\kappa)$ one has $M^*_+(\alpha)  M_+(\kappa) h=M_+(\overline{\alpha})P_{H_+}(\kappa h)=P_{H_+}M(\overline{\alpha})P^*_{H_+} P_{H_+}(\kappa h)=P_{H_+}(\overline{\alpha}\kappa h)-P_{H_+}\big(\overline{\alpha}P^*_{H_-}P_{H_-}(\kappa h)\big)=M_+(\overline{\alpha}\kappa)h-0$, since $\overline{\alpha}P^*_{H_-}P_{H_-}(\kappa h)\in H_-$. --- For the translational invariance check $S_b=
W(\e^{\i b(\cdot)})$, $\dom W(\e^{-\i b(\cdot)}\kappa)=\dom W_\kappa$ and apply the foregoing results.\\
(e) Here we prove the case $p=2$ and the implication \,$\Leftarrow$\, for $p=\infty$ and $p=1$. The remainder is shown in the proof of (f).\\
\hspace*{6mm}
Let $h\in\dom M_+(\kappa)\setminus \{0\}$. Then $\kappa h\in L^2(\R)$, whence $\kappa \e^{\ln|h|}\in L^2(\R)$. By   (\ref{MHF}), $\ln|h|/(1+x^2)$ is integrable. For the converse implication  
put $j':=1_{\{j\le 0\}} j- |x|^{1/2}$. Then  $j'/(1+x^2)$ is integrable, $\e^{j'}$ is square-integrable. By  (\ref{MHF}) there is $h\in H_+$ with $|h|=\e^{j'}$. Moreover, $\kappa \e^{j'}$ is square-integrable since $j'\le j$. So $h\in\dom M_+(\kappa)\setminus \{0\}$.   --- Next turn to \,$\Leftarrow$\, for $p=\infty$. For $j'$ from above    $|\kappa| \e^{j'}\le |\kappa| \e^j \e^{-|x|^{1/2}}$ is square-integrable. Hence $\kappa$ is proper by the case $p=2$  just shown. --- Now consider  \,$\Leftarrow$\, for $p=1$. Due to the assumption  
$|\kappa|\e^j(1+x^2)^{-1}$ is integrable, whence $j'(1+x^2)^{-1}$ is integrable for 
$j':=\ln(1+|\kappa|\e^j)$. Therefore $(1+|\kappa|\e^j)\e^{-j'}=1$ implying $|\kappa|\e^{j-j'}\le1$, whence the result by  \,$\Leftarrow$\, for $p=\infty$.\\
(f) Let $\kappa^2$ be proper.  Then $1+|\kappa|^2$ is proper and $|\kappa|\le1+|\kappa|^2$, whence $\kappa$ is proper. --- Now let $\kappa$ be proper.  Then by (e)($p=2$), $\kappa \e^j\in L^2(\R)$ for some  real-valued $j$ with integrable $j/(1+x^2)$. Since $\kappa':=|\kappa|^2\e^{2j}$ is integrable, $\kappa'$ is proper  by (e)($\Leftarrow$\, for $p=1$). Then (e)($p=2$) yields $\kappa' \e^{j'}\in L^2(\R)$ for some real-valued $j'$ with integrable $j'/(1+x^2)$. Hence $|\kappa|^2 \e^{2j+j'}\in L^2(\R)$, whence the claim by (e)$(p=2)$.\\
\hspace*{6mm}
The general case is easily reduced to the claim that, for $r>1$ and $\kappa\ge 0$, $\kappa$ is proper if and only if $\kappa^r$ is proper.  So let $\kappa$ be proper.  Let $n\in\N$ satisfy $r\le 2^n$. By the foregoing result $\kappa^{2^n}$ is proper. Then $1+\kappa^{2^n}$ is proper and $\kappa^r\le 1+\kappa^{2^n}$. Hence $\kappa^r$  is proper. Conversely, if $\kappa^r$ is proper, then $1+\kappa^r$ is proper and $\kappa\le 1+\kappa^r$, whence $\kappa$ is proper.\\
\hspace*{6mm}
Now we complete the proof of (e). Consider first the case $p\in]0,\infty[$.
Let $\kappa \e^j\in L^p(\R)$ for some $j$. Then 
$|\kappa|^{p/2} \e^{\frac{p}{2}j}\in L^2(\R)$. Hence $|\kappa|^{p/2}$ is proper as shown in (e)($p=2$). The foregoing result applies, whence $\kappa$ is proper. The converse follows in the same way due to  (e)($\Rightarrow$\, for $p=2$). --- Now let $p=\infty$ and let $\kappa$ be proper. As just shown, $\kappa e^j$ is integrable for some real-valued $j$ with integrable $j(1+x^2)^{-1}$. Then $|\kappa|\e^{j''}\le 1$ for $j'':=j-\ln(1+|\kappa|\e^j)$ as shown in the proof of (e)($\Leftarrow$\, for $p=1$).\\
\hspace*{6mm}
Finally, let $\kappa_1$, $\kappa_2$ be proper. Apply (e) for $p=\infty$. So $|\kappa_i|\le\e^{j_i}$, $i=1,2$. Then
$|\kappa_1\kappa_2|\le\e^{j_1+j_2}$. Assume without restriction $j_i\ge 0$, $i=1,2$. Then $|\kappa_1+\kappa_2|\le 2\e^{j_1+j_2}$. \\ 
(g) Let $\kappa$ be proper. Then $1+|\kappa|$ is proper, and by (e)($p=\infty$) one has $1+|\kappa|\le \e^j$  for some real-valued $j$ with integrable $\frac{j}{1+x^2}$.
Then $\ln(1+|\kappa|)\le j$, whence the claim. Conversely let $\frac{j}{1+x^2}$ be integrable for $j:=\ln(1+|\kappa|)$. Then $1+|\kappa|=\e^j$, whence $|\kappa|\le \e^j$ and $\kappa$ is proper by (e)($p=\infty$).\\
(h) $\ln(1+|\kappa(x)|)\ge\ln|\kappa(x)|\ge \ln(\delta)+\lambda |x|$ for all $x\ge a$ or $x\le -a$. Hence $\frac{\ln(1+|\kappa|)}{1+x^2}$ is not integrable. Apply (g).\\
(i) Let $h\in\dom M_+(\kappa)\setminus \{0\}$. Then $\kappa h\in L^2(\R)$, whence $\kappa \e^j\in L^2(\R)$ for $j:=\ln|h|$. By   (\ref{MHF}), $j/(1+x^2)$ is integrable.\\ 
\hspace*{6mm}
Recall (\ref{RTO}). Put $c:= C^{-1}$.
Check that $j\circ c\ge 0$ is integrable on the torus $\mathbb{T}$. Hence $U(w):=\exp\big(\frac{1}{2\pi}\int_0^{2\pi}\frac{\e^{\i t}+w}{\e^{\i t}-w}\,j\circ c(\e^{\i t})\d t\big)$ is an outer function on  the disc $\mathbb{D}$ (see e.g.\,\cite[Chapter 3]{A14}). 
So $U$ converges non-tangentially a.e. to a function $u$ on  $\mathbb{T}$ satisfying $|u|=\e^{j\circ c}\in L^2(\mathbb{T})$, whence  $U \in  H^2(\mathbb{D})$. Therefore, as known (see e.g. \cite[sec.\,3]{V03}), $U H^\infty(\mathbb{D})$ is dense in $H^2(\mathbb{D})$. This implies that $uH^\infty(\mathbb{T})$ is dense in $H^2(\mathbb{T})$.\\
\hspace*{6mm} 
Put  $h':=u\circ C$. Note  $|h'|=\e^j=|h|$.
The above  result is transferred to $H_+$ by   
$\Gamma$ in (\ref{RTO}).
Accordingly, $\frac{1}{x+\i}h'\in H_+$ and $\frac{1}{x+\i}h' H_+^\infty$ is dense in $H_+$. The latter is contained in $\dom M_+(\kappa)$ by (c), since $\kappa \frac{1}{x+\i}h'$ is square-integrable and hence  $\frac{1}{x+\i}h'\in\dom M_+(\kappa)$. Finally $h'\in H_+$ by (\ref{AL1})(c).\\
(j) Note $0<\rho\le 1$. So $h\in\dom M_+(\kappa)=\dom M_+(\rho^{-1})$ $\Leftrightarrow$ $\exists f\in L^2(\R)$ such that $h=\rho f\in H_+$. The latter means $0=\langle k,\rho f\rangle=\langle \rho k,f\rangle$ $\forall$ $k \in H_-$. Hence $P^*_{H_+}\dom M_+(\kappa)=M(\rho)\big(M(\rho)H_-\big)^\perp$. --- Now let $h_0\in H_+$ with $h_0\perp M(\rho)\big(M(\rho)H_-\big)^\perp$. The latter means $0=\langle h_0,\rho f\rangle=\langle\rho h_0,f\rangle$ $\forall$ $f\in \big(M(\rho)H_-\big)^\perp$ and hence equivalently $\rho h_0 \in \overline{M(\rho)H_-}$.  --- Since by (i) $\dom M_+(\kappa)$ is either dense if $\kappa$ is proper or trivial if $\kappa$ is not proper, the result follows.\\
(k) Cf.\;the proof for bounded $\kappa$ in \cite[2.8]{D79}. Keep (\ref{HS}) in mind. Since $W_\kappa-\lambda I=W(\kappa-\lambda1_\R)$ assume without restriction $\lambda=0$. Suppose $W_\kappa u=0$, $W_{\overline{\kappa}}v=0$. Set $h_+:=\F^{-1}P_+^*u$. Then $h_+\in H_+$ and $h_-:=\kappa h_+\in H_-$. Set $k_-:=\overline{\F^{-1}P_+^*v}$. Then $k_-\in H_-$ and $\overline{\kappa}\overline{k}_-\in H_-$, whence $k_+:=\kappa k_-\in H_+$.\\
\hspace*{6mm}
 Note that $j:=h_-k_-=\kappa h_+k_-=h_+k_+$.  Hence there is a holomorphic $\chi:\C\setminus\R\to\C$    such that its partial maps satisfy $\chi_y\in L^1(\R)$, $K:=\sup\{\norm{\chi_y}_1:y\ne 0\}<\infty$, $\chi_y\to j$ pointwise a.e. and $\norm{\chi_y-j}_1\to 0$ for $y\to 0$.
By a standard argument (see also  \cite[Theorem II]{C44}) $\chi$ extends to an entire function still called $\chi$ with $\chi|_\R=j$ a.e. Fix $z\in\C$, $|z|>1$. We use the representation  $\pi \chi(z)=\int_D \chi(z+w)\d^2w$ were $D$ denotes the disc with center $0$ and radius $1$. 
Then  $\pi|\chi(z)|\le \int_{-1}^1\int_{-\infty}^\infty |\chi\big(x+u+\i(y+v)\big)|\d u\d v\le 2K$ so that $\chi$ is constant equal to $0$, whence $h_-k_-=0$, $h_+k_+=0$. (An alternative argument implying this result uses \cite[Theorem 76]{T48}, by which $(2\pi)^{1/2}\F_{L^1}(h_\pm k_\pm)=\F h_\pm\star \F k_\pm$. Accordingly, $\F h_-\star \F k_-=\F h_+\star \F k_+=0$, since $\F h_-$,  $\F k_-$ vanish on $[0,\infty[$ and $\F h_+$,  $\F k_+$ vanish on $]-\infty,0]$.) The cases $k_-=0$ or  $h_+=0$ are trivial. Otherwise  $h_+\ne 0$ and $\kappa h_+=h_-=0$, whence $\kappa=0$ a.e. \\
\hspace*{6mm}
Now let $\kappa$ be real, let $\lambda\in\C$, and let $g\in \dom W_\kappa$ satisfy $W_\kappa g=\lambda g$. Then $\lambda\,\langle g,g\rangle =\langle g,W_\kappa g\rangle = \langle W_\kappa g, g\rangle=\overline{\lambda}\,\langle g,g\rangle$. Hence $g=0$, since otherwise $\lambda\in\R$ would contradict the foregoing result.\\
(l) Since $W_\kappa-\lambda I=W(\kappa-\lambda1_\R)$ and  $\kappa^{-1}(\{\lambda\})=(\kappa-\lambda 1_\R)^{-1}(\{0\})$     assume without restriction    $\lambda=0$. Put $E:=\kappa^{-1}(\{0\})$. Suppose $W_\kappa g=0$. Then $h_+:=\mathcal{F}^{-1}P_+^*g \in H_+$ and $h_-:=\kappa h_+\in H_-$. Since $h_-$ vanishes on the non-null set $ E$,  $h_-=0$ follows. Hence $h_+$ vanishes on the non-null set $\R\setminus E$ implying $h_+=0$ and hence $g=0$.\\
(m) Put $\beta:=\kappa_1-  \kappa_2$. Then $W_\beta|_D =0$ for some dense $D\subset L^2(\R_+)$. Let $k\in\dom W_\beta^*$. Then for all $g\in D$ one has $0=\langle k,W_\beta g\rangle=\langle W_\beta^*k,g\rangle$, whence $W_\beta^*k=0$. Since $W_{\overline{\beta}}\subset W_\beta^*$, it follows $W_{\overline{\beta}}|_D=0$. So $\beta$ is almost constant by (k), whence $\beta=0$ a.e.\\
  (n) Suppose  $W_\kappa\subset W_\kappa^*$.  Since generally $W_{\overline{\kappa}}\subset W_\kappa^*$ and $\dom W_\kappa=\dom W_{\overline{\kappa}}$, it follows $ W_\kappa= W_{\overline{\kappa}}$, whence $\kappa=\overline{\kappa}$ a.e. by (m).  The converse is  obvious.\\
(o) By (n) $\kappa$ is real. First suppose that $\kappa$ is bounded below with maximal lower bound $a$. Put $\varkappa:=\kappa-a\ge 0$. Then 
$\langle h, M_+(\varkappa)h\rangle=\int \varkappa |h|^2 \d x\ge 0$, 
 whence $\langle h, M_+(\kappa)h\rangle\ge a\norm{h}^2$ $\forall$ $h\in\dom M_+(\varkappa)=\dom M_+(\kappa)$. So $a$ is a lower bound for $M_+(\kappa)$.\\
\hspace*{6mm}
Now let $M_+(\kappa)$ be bounded below with maximal lower bound $\alpha$. Put $\varkappa:=\kappa-\alpha$. Then 
$\langle h, M_+(\varkappa)h\rangle\ge 0$ $\forall$ $h\in\dom M_+(\varkappa)=\dom M_+(\kappa)$, and 
for  $\epsilon>0$ there is $h_1\in\dom M_+(\varkappa)$, $\norm{h_1}=1$  with $\langle h_1, M_+(\varkappa)h_1\rangle< \epsilon $. Put $A:=\{x\in\R:\varkappa(x)\ge 0\}$. Hence $\int1_A\varkappa\, |h_1|^2 \d x<\epsilon$. Assume that  $B:=\R\setminus A$ is not null set. Let $h_2\in\dom M_+(1_B\varkappa)\supset \dom M_+(\varkappa)$, $h_2\ne 0$. Then 
$c:=\int 1_B\varkappa\,|h_2|^2 \d x<0$. Choose $\epsilon<|c|$. The function $f:=1_A|h_1|+1_B|h_2|$ satisfies $|\operatorname{ln}f(x)|\le |\operatorname{ln}|h_1(x)|\,|+|\operatorname{ln}|h_2(x)|\,|$. Since $h_1,h_2\in H_+$ it follows that $\operatorname{ln}f(x)/(1+x^2)$ is integrable, whence there is $h\in H_+$ with $|h|=f$. Note  that $h\in\dom M_+(\varkappa)$, since $|\varkappa h|=|1_A\varkappa \,h_1|+|1_B\varkappa \,h_2| \in L^2(\R)$. It follows the contradiction $\langle h, M_+(\varkappa) h\rangle=\int \varkappa |h|^2\,\d x\le c+\epsilon<0$. Therefore $B$ is a null set,  whence $\varkappa\ge 0$ a.e. and  $\alpha$ is an essential lower bound of $\kappa$. This proves $a=\alpha$.\\
\hspace*{6mm}
$M_+(\kappa)$ not being bounded above by assumption, the numerical range $R$ of $M_+(\kappa)$ is not bounded above. As $R$ is convex one  infers $]\alpha,\infty[\subset R\subset [\alpha,\infty[$.
It remains to show $\alpha\not\in R$ or, equivalently,  that $\langle h,M_+(\varkappa)h\rangle=0$ implies $h=0$, where $\varkappa\ge 0$ is not almost zero. Indeed, $0=\langle h,M_+(\varkappa)h\rangle=\int \varkappa |h|^2\d x$ $\Rightarrow$ $h$ vanishes on a non-null set $\Rightarrow$ $h=0$.\\
(p) follows readily from (o).\qed\\

 In view of (\ref{DBWHO})(k) we recall that
 self-adjoint bounded WH, which are not a multiple of $I$, i.e., $\kappa$  real bounded not almost constant,  are  even absolutely continuous. Indeed, these operators are  Hilbert space isomorphic to self-adjoint bounded Toeplitz operators (see (\ref{RTO}), or e.g.  \cite[9.5(e)]{BS90}, \cite[3.3.2\,(13)]{C19}), which by Rosenblum \cite{R65} are absolutely continuous. Actually, as already mentioned,  it follows from \cite{R65} that for   real bounded below not almost constant $\kappa$,  for which $(1+x^2)^{-1}\kappa$ is integrable, the Friedrichs extension $\tilde{W}_\kappa$ (\ref{SBS}) of $W_\kappa$    is   absolutely continuous.\\
\hspace*{6mm}
Section \ref{SBWHO} is concerned with the case that 
 $\kappa$ is proper real and semibounded. This is the general case that $W_\kappa$ is  densely defined symmetric  semibounded (\ref{DBWHO})(n),(o). The natural self-adjoint extension  $\tilde{W}_\kappa$ is studied in  (\ref{RSBWHO}), (\ref{SBS}). \\
\hspace*{6mm}
 If $\kappa$ is proper real and even (i.e. $\kappa(-x)=\kappa(x)$) then $W_\kappa$ is densely defined symmetric and  has a self-adjoint extension. This holds true  since $L^2(\R_+)\to L^2(\R_+)$, $g\mapsto \overline{g}$ is a conjugation, which  leaves  $\dom W_\kappa$ invariant and satisfies 
$W_\kappa \overline{g}=\overline{W_\kappa g}$ (see  \cite[Theorem X.3]{RS75}). If $\kappa$ is odd instead of even then in general $W_\kappa$ has no self-adjoint extension. Examples are furnished by  real rational symbols  as e.g. $\kappa(x)=x$. In (\ref{DIWHO}) an explicit description of the deficiency spaces of $W_\kappa$ for real rational $\kappa$ are given yielding further examples of densely defined symmetric WH with self-adjoint extensions.\\

Concluding this section we deal with the unilateral translation invariance (\ref{DBWHO})(d) of WH.  We are inspired by \cite[sec.\,2.3)]{Y17} which treats the bounded case (\ref{ULTICC}). Observe   the easily verifiable relation 
\begin{equation}\label{IUT} 
T_b^*P_+^*W_\kappa P_+T_b= M(1_{]-b,\infty[})\F M(\kappa)\F^{-1}M(1_{]-b,\infty[})
\end{equation}
 where $T_b$, $b\in\R$ denotes the unitary one-parameter group of translations $T_bf(x):=f(x-b)$ on $L^2(\R)$.  It shows again the invariance  
\begin{equation}\label{ULTI}
W_\kappa=S^*_b W_\kappa S_b
\end{equation}
 under the unilateral translations $S_b$, $b\ge 0$ (see  (\ref{DBWHO})(d)), since $S_b=P_+T_bP_+^*$. Recall 
 $S_b=W(\e^{\i b(\cdot)})$. 
 Moreover it implies that $T_b^*P_+^*W_\kappa P_+T_bf$ converges as $b\to\infty$ if $f\in T_a^*P_+^*(\dom W_\kappa)$  for some $a\ge 0$, yielding the limit $\F M(\kappa)\F^{-1}f$.  In particular $\{\norm{W_\kappa S_bg}:b\ge0\}$ is bounded for every $g\in\dom W_\kappa$. If $W_\kappa$ is densely defined then also $\{\norm{W^*_\kappa S_bg}:b\ge0\}$ is bounded
 as $ W^*_\kappa|_{\dom W_\kappa}=W_{\overline{\kappa}}$.
 
\begin{The}\label{ULTIC} Let $A$ be a densely defined operator in $L^2(\R_+)$ satisfying  \begin{equation}\label{EULTIC}
 A\subset S_b^*AS_b\quad \forall \;b\ge0
 \end{equation}
 If $\dom A\subset
 \dom A^*$ and 
$\{\norm{AS_bg}:b\ge 0\}$ and $\{\norm{A^*S_bg}:b\ge 0\}$ are bounded for every $g\in\dom A$ then there is a WH $W_\kappa$ extending $A$.
 \end{The}\\ 
{\it Proof.} 
(i) For $a\ge 0$ put $D_a:=T_a^*P_+^*(\dom A)$ and  $D:=\bigcup_{a\ge0} D_a$. Then  $D$ is a dense translation invariant subspace of  $L^2(\R)$.\\
 \hspace*{6mm} Indeed, let $f\in L^2(\R)$ and $\varepsilon>0$. Since $f_a:=1_{]-a,\infty[}f \to f$ for $a\to \infty$ in the mean, one has $\norm{f-f_a}\le\varepsilon/2$ for some $a\ge0$. As $T_af_a=P_+^*P_+T_af$, there is $g\in \dom A$ with $\norm{P_+T_af_a-g}\le\varepsilon/2$. Hence $\norm{f_a-T_a^*P_+^*g}\le\varepsilon/2$.  It follows that $D$ is dense.\\
\hspace*{6mm} 
Note  that $S_b(\dom A)\subset \dom A$ for $b\ge 0$ due to $A\subset S_b^*AS_b$. Then $D_a\subset D_{b}$ for $a\le b$. Indeed, put $c:=b-a\ge0$. Then  check $T_a^*P_+^*=T^*_{b}P^*_+S_c$, whence $D_a= T^*_{b}P^*_+S_c(\dom A)\subset T^*_{b}P^*_+(\dom A)=D_{b}$. It follows that $D$ is a  subspace of $L^2(\R)$.\\
 \hspace*{6mm}
 For the translation invariance of $D$  it suffices to show $T_c^*P_+^*g\in D$ for $c\in\R$, $g\in \dom A$. If $c\ge 0$ this is obvious. Let $c<0$. Then $T_c^*P_+^*g=T_{-c}P_+^*g=P_+^*P_+T_{-c}P_+^*g=P_+^*S_{-c}g\in D$ as $S_{-c}g\in \dom A$.\\
 \hspace*{6mm}
(ii) Let $f\in D$. Then $f=T_a^*P_+^*g$ for some  $a\ge 0$, $g\in \dom A$, whence $f_b:=T_b^*P_+^*A P_+T_bf=T_b^*P_+^*A S_{b-a}g$ is well-defined for all $b\ge a$.  We are going to show that $\lim_{b\to \infty}f_b$ exists.\\
\hspace*{6mm}
 Indeed, let $c\ge b\ge a$. Then 
$\norm{f_c-f_b}^2=\norm{f_c}^2+\norm{f_b}^2-\langle f_c,f_b\rangle-\langle f_b,f_c\rangle$ with $\norm{f_c}^2=\norm{AP_+T_cf}^2$, $\norm{f_b}^2=\norm{AP_+T_bf}^2$, $\langle f_c,f_b\rangle=
\langle S^*_{c-b}AP_+T_{c-b}T_bf,AP_+T_bf\rangle=\langle S^*_{c-b}AP_+T_{c-b}(P_+^*P_+T_bf),AP_+T_bf\rangle=
\langle S^*_{c-b}AS_{c-b}P_+T_bf,AP_+T_bf\rangle=||AP_+T_bf||^2$ by (\ref{EULTIC}) since $P_+T_bf\in\dom A$. Hence $\norm{f_c-f_b}^2=\norm{AP_+T_cf}^2-\norm{AP_+T_bf}^2$. This implies that $b\mapsto \norm{AP_+T_bf}^2$ is increasing. Being bounded by the assumption it follows $\norm{f_c-f_b}^2\to 0$ for $b,c\to \infty$ so that $\lim_{b\to \infty}f_b$ exists.
\\
 \hspace*{6mm}
(iii) Thus $C_1f:=\lim_{b\to \infty}f_b$ defines an operator on $D$. It is translation invariant since $T_c^*C_1T_cf=\lim_{b\to\infty}T_c^*T_b^*P_+^*A P_+T_bT_cf=\lim_{b\to\infty}T_{c+b}^*P_+^*A P_+T_{c+b}f=\lim_{b'\to\infty}T_{b'}^*P_+^*A P_+T_{b'}f=C_1f$.\\
 \hspace*{6mm} 
(iv) Also $C_2f:= \lim_{b\to\infty}T_b^*P_+^*A^* P_+T_bf$ exists for $f\in D$ thus defining an operator $C_2$ on $D$. This result follows replacing $A$ in (ii) by $A^\#:=A^*|_{\dom A}$. It remains to verify $A^\#\subset S_b^*A^\#S_b$, $b\ge0$. Indeed, for $g,g'\in \dom A$ one has $\langle g,S_b^*A^\#S_bg'\rangle=\langle S_bg,A^*S_bg'\rangle=\langle AS_bg,S_bg'\rangle$ since $S_bg\in\dom A\subset\dom\overline{A}=\dom A^{**}$. Hence $\langle g,S_b^*A^\#S_bg'\rangle=\langle S_b^*AS_bg,g'\rangle=\langle Ag,g'\rangle=\langle g,A^\#g'\rangle$ by (\ref{EULTIC}) and $g,g'\in \dom A$. This implies $S_b^*A^\#S_bg'=A^\#g'$, whence the claim.\\
\hspace*{6mm}
(v) Obviously $C_2\subset C_1^*$. Hence $C_1^*$ is densely defined and the closure $C:=\overline{C}_1$ exists. Clearly translation invariance $C=T_b^*CT_b$  holds. Equivalently  $\F^{-1} C\F$ commutes with $M(\e^{\i b(\cdot)})=\F^{-1} T_b\F$ for all  $b\in \R$. Thus $\F^{-1} C\F=M(\kappa)$ for some measurable function $\kappa$. Hence $P_+CP_+^*=W_\kappa$. Finally, for $g\in\dom A$ one has $P_+^*g\in D$ 
and $W_\kappa g=P_+C_1P_+^*g=\lim_{b\to\infty}P_+T_b^*P_+^*AP_+T_bP_+^*g=\lim_{b\to\infty}S_b^*AS_bg
=\lim_{b\to\infty}Ag=Ag$.\qed

\begin{Cor}\label{ULTICC}
 Let $A$ be a bounded operator on $L^2(\R_+)$.  Then  $A$ is a WH if and only if $$A = S_b^*AS_b\quad \forall \;b\ge0$$ 
 \end{Cor}
{\it Proof.} It remains to observe that $A\subset W_\kappa$ by (\ref{ULTIC}) implies $A=W_\kappa$.\qed\\
 For (\ref{ULTICC}) see also  \cite[(2.10)]{Y17}, where  
 the existence of $\lim_{b\to \infty}T_b^*P_+^*A P_+T_bf$ (see (ii) of the proof of  (\ref{ULTIC})) is not proven.

\section{Rational Symbols}\label{RS}

WH for rational symbols $\kappa=\frac{P}{Q}\big|_\R$ with polynomials $P\ne 0$, $Q\ne 0$ permit some more  general analysis.  According to (\ref{DBWHO})(b) they are densely defined. In (\ref{RSKCK})  we show that they are  closed and  we determine their domains, ranges, and kernels and deficiency spaces, which are finite dimensional, and their spectral and Fredholm points. In particular, in the symmetric case, i.e., for a real rational symbol the deficiency spaces and indices   are explicitly available (\ref{DIWHO}).\\ 
\hspace*{6mm}
Mostly we will omit $|_\R$ indicating the restriction on $\R$. A polynomial with a negative degree is the null function.

\begin{Lem}\label{AL1} Let $P\ne 0$ and $Q\ne 0$ be  polynomials.

\emph{(a)} Let $P$ and $Q$ have no common zeros. Then $\frac{P}{Q}\in H_+$ ($\in H_-$) if and only if $\deg P<\deg Q$ and  all zeros of  $Q$ are in the lower  (upper) half-plane.

\emph{(b)}  Let $h\in H_+\setminus\{0\}$ such that $\frac{P}{Q}h\in H_-$. Then there is a polynomial $R$ with $\operatorname{deg} R<\min\{\deg P, \deg Q\}$ such that all zeros of $P$ in the closed upper half-plane   as well as all zeros of $Q$ in the closed lower half-plane  are zeros of $R$ and such that $h=\frac{R}{P}$. Conversely, it is obvious that $h_+:=\frac{R}{P}$ and  $h_-:=\frac{R}{Q}$ satisfy $h_\pm\in H_\pm$ and $\frac{P}{Q} \,h_+=h_-$.

\emph{(c)} Let $Q$ have no zeros in the upper  half-plane. If $h\in H_+$ and $\frac{P}{Q}h\in L^2(\R)$, then $\frac{P}{Q}h\in H_+$.

\emph{(d)} Let $P$ and $Q$ have no common zeros. 
 Suppose that $h\in H_+$ and $\frac{P}{Q}h\in H_+$. Then $h/Q\in H_+$.
 \end{Lem}\\
{\it Proof.} (i) Regarding (a) put $h:=\frac{P}{Q}\big|_\R$. Obviously $h\in L^2(\R)$ if and only if $\deg P<\deg Q$ and  $Q$ has no real  zeros. Moreover, if all zeros of $Q$ are in the lower half-plane, then $\psi:=\frac{P}{Q}$ is holomorphic in the upper half-plane with bounded $\{\norm{\psi_y}_2: y>0\}$, whence $h\in H_+$. Conversely, let $\phi$ be the holomorphic function on the upper half-plane associated with $h\in H_+$. Then 
$\phi$ converges to $h$ non-tangentially a.e.,
and $\psi$ is meromorphic without real poles, whence  $\psi(z)\to h(x)$ for $z\to x\in\R$.
Hence $\phi$ and $\psi$ coincide on the upper half-plane by \cite[IV 2.5]{Pr65} so that $\psi$ has no poles there.

(ii) As to (b) we prove $h=R/P$ for some polynomial $R$ supposing that $Q$ has no real zeros and $\deg Q\ge \deg P$.  \\
\hspace*{6mm}
 Let $\phi$ and  $\psi$ be the holomorphic functions on the upper half-plane and lower half-plane associated with $h$ and $\frac{P}{Q}h$, respectively. There is $\delta>0$ such that $B:=P/Q$ is holomorphic and bounded in the strip  $\{z:-\delta<\operatorname{Im} z<\delta\}$, and hence $B\phi$ is holomorphic on  $\{0<\operatorname{Im} z<\delta\}$. 
Since $B$ is bounded, one  easily infers $\int_J|(B\phi)(x+\i y)-\psi(x-\i y)|\d x\to 0$, $0<y\to 0$ for any bounded interval $J$. Then by \cite[Theorem II]{C44} there is a holomorphic function $\chi$  on $\{z:\operatorname{Im} z<\delta\}$ extending $B\phi$  and $\psi$.  So $Q\chi$ is still holomorphic on $\{z:\operatorname{Im} z<\delta\}$ coinciding  with $P\phi$ on $\{z:0<\operatorname{Im} z<\delta\}$ and with $Q\psi$  on  the lower half-plane. Hence there is an entire function $R$ extending $P\phi$ on the upper and $Q\psi$ on the lower half-plane. Introduce $S$ being equal to $P$ on the upper half-plane and equal to $Q$ on the lower half-plane. Analogously define  $\varSigma$ with respect to $\phi$ and $\psi$.\\
\hspace*{6mm}
Fix $z\in\C$, $|z|>1$. We use the representation  $R(z)=\frac{1}{\pi}\int_D R(z+w)\d^2w$ were $D$ denotes the disc with center $0$ and radius $1$. 
Then  
$$|R(z)|^2\le \pi^{-2}\int_D|S(z+w)|^2\d^2w\;\int_{-1}^1\int_{-\infty}^\infty |\varSigma\big(x+u+\i(y+v)\big)|^2\d u\d v$$
The first integral is easily estimated $\le$ constant $|z|^{2n}$ with $n:= \deg Q$. The double integral is bounded independently of $z$, since $\int |\varSigma(u+\i v)|^2\d u\le K$  for all  $v\ne0$ with some constant $K<\infty$. Therefore $R$ is a polynomial  and $\phi=\frac{R}{P}$.

(iii) Next we show  $Ph\in H_+$ for  $h\in H_+$  if $Ph\in L^2(\R)$.\\
\hspace*{6mm}
Let $a\in\C$ be a zero of $P$ and write $P=(x-a)P'$. Note that $(x-a)h\in L^2(\R)$ since $|x-a|\le c|P(x)|$,  $x\in\R\setminus J$ for some bounded interval $J$ and constant $c$. Hence it suffices to prove $(x-a)h\in H_+$ and proceed with $P'$ in place of $P$. Actually the claim is $xh\in H_+$ and hence equivalently $(x+\i)h\in H_+$.\\
\hspace*{6mm}
Let $k\in H_-$. Then $\frac{k}{x-\i}\in H_-$ and $\langle \frac{k}{x-\i},(x+\i)h\rangle=\langle k,h\rangle=0$. To conclude the proof obviously it suffices to show that $\{\frac{k}{x-\i}:k\in H_-\}$ is dense in $H_-$. So let $k_0\in H_-$ satisfy $\langle k_0,\frac{k}{x-\i}\rangle=0$ for all $k\in H_-$.
Then $\langle \frac{k_0}{x+\i},k\rangle=0$ implying $\frac{k_0}{x+\i}\in H_+$ or equivalently $\frac{\overline{k}_0}{x-\i}\in H_-$, whence $k_0=0$ by (ii).

(iv) Let $Q$ have no zeros in the upper  half-plane. The claim is $\frac{1}{Q}h\in H_+$ for $h\in H_+$ if $\frac{1}{Q}h\in L^2(\R)$.\\
\hspace*{6mm}
Write $Q=q\,Q_<$, where the zeros of $q$ are exactly the real zeros of $Q$. Let $q_\epsilon(z):=q(z+\i\epsilon)$ for $\epsilon>0$ and put $Q_\epsilon:=q_\epsilon Q_<$. Note $|q(x)/q_\epsilon(x)|<1$, $x\in\R$. Hence $h/Q_\epsilon\in L^2(\R)$. Since $1/Q_\epsilon$ is bounded on the upper half-plane, $h/Q_\epsilon\in H_+$. Moreover, $h/Q_\epsilon\to h/Q$ for $\epsilon\to 0$ pointwise and in the mean, whence $h/Q\in H_+$.

(v) It follows the proof of (c). 
Without restriction let $P$ and $Q$ be without common zeros. Then  $\frac{P}{Q}h\in L^2(\R)$ implies  $\frac{1}{Q}h\in L^2(\R)$. Indeed, let $K$ be a compact neighborhood of the real zeros of $Q$ containing no real zero of $P$. Then $|h/Q|$ is bounded by  $C|P\,h/Q|$ on K and by $c|h|$ on $\R\setminus K$ for some finite constants $C,c$, whence $\frac{1}{Q}h\in L^2(\R)$. --- 
Now $\frac{1}{Q}h\in H_+$ by (iv) and hence $P\frac{1}{Q}h\in H_+$ by (iii).

(vi) We proceed with the proof of (b). 
Without restriction let $P$ and $Q$ be without common zeros. Write $Q=q\,Q_0$,  where the zeros of $q$ are exactly the real zeros of $Q$. Then  $\frac{P}{Q}h\in L^2(\R)$ implies  $\frac{1}{q}h\in L^2(\R)$, cf.\,(v). Hence $h':=h/q\in H_+$ by (iv). Next let $P=p\,P_0$ with $\deg P_0=\deg Q_0$ if $\deg P>\deg Q_0$ and  $p=1$ otherwise. Then $h'':=ph'\in L^2(\R)$, since $h''=\frac{Q_0}{P_0} \frac{P}{Q} h$, where in the case $p\ne 1$ the factor  $\frac{Q_0}{P_0}$ is bounded outside a bounded interval. By (iii) this implies $h''\in H_+$. By assumption $\frac{P_0}{Q_0}h''\in H_-$. Therefore $h''=R_0/P_0$ by (ii) for some polynomial $R_0$. It follows $h=R/P\in H_+$ for $R:=qR_0$ and hence $R/Q\in H_-$. The proof is accomplished applying (a) proved in (i).

(vii) As to the proof of (d)  assume first that $Q$ has no real zeros. Let $\phi$ and  $\psi$ be the holomorphic functions on the upper half-plane  related to $h$ and $\frac{P}{Q}h$, respectively. Then $\frac{P}{Q}\phi$ is meromorphic on the upper half-plane and converges non-tangentially to  $ \frac{P}{Q}h$ a.e. 
Since $\psi$ does the same,
according to  \cite[IV 2.5]{Pr65}, $\psi=\frac{P}{Q}\phi$ holds. Hence $\phi/Q$ is holomorphic on the upper half-plane 
with $(\phi/Q)_y\to h/Q$ for $0<y\to 0$ a.e. 
Let $0<\delta<c$ such that $C:=[-c,c]\times \i[\delta, c]$ is a neighborhood  of the zeros of $Q$ in the upper half-plane. 
Then $|1/Q|$ is bounded by some constant $L$ on $\{z:\operatorname{Im}z\ge0\}\setminus C$,  and   $|\phi/Q|$ is bounded on $C$ by some   $M$. Recall that $\norm{\phi_y}^2_2$ is bounded for $y>0$ by some $K$. Then 
$\int|(\phi/Q)(x+\i y)|^2\d x\le L^2\int |\phi(x+\i y)|^2\d x+\int_{-c}^cM^2\d x\le L^2K+2cM^2$ for all $y>0$.
Finally, for $0<y<\delta$, $\int|(\phi/Q)(x+\i y)-(h/Q)(x)|^2\d x\le L^2\int|\psi(x+\i y)-h)(x)|^2\d x\to 0$ for $y\to 0$, whence $h/Q\in H_+$.\\
\hspace*{6mm}
Now we turn to the general case. Note  $h/Q\in L^2(\R)$, cf.\,(v). Write $Q=q\,Q_0$, where the zeros of $q$ are exactly the real zeros of $Q$. Let $q_\epsilon(z):=q(z+\i\epsilon)$  for $\epsilon>0$ and put $Q_\epsilon:=q_\epsilon Q_0$. Note $|q(z)/q_\epsilon(z)|<1$ on the upper half-plane. Therefore $\frac{P}{Q_\epsilon}h=\frac{q}{q_\epsilon}\frac{P}{Q}h\in H_+$, where $Q_\epsilon$ has no real zeros. Moreover, for $\epsilon>0$ small enough, $P$ and $Q_\epsilon$ have no common zeros. Hence the foregoing result applies so that 
$h/Q_\epsilon\in H_+$. Now $h/Q_\epsilon=\frac{q}{q_\epsilon}h/Q\to h/Q$ for $\epsilon\to 0$ in the mean implying $h/Q\in H_+$. 
\qed\\

Recall that a densely defined closed operator between Banach spaces with finite dimensional kernel  and cokernel is called a Fredholm operator if its range is closed  (cf.\,\cite{S67}).

\begin{The}\label{RSKCK} Let $\kappa=\frac{P}{Q}$ be a rational function, where the polynomials $P$ and $Q$ have no common zeros.  Then $M_+(\kappa)$ is densely defined and closed and
\begin{itemize}
\item[\emph{(a)}] $\dom  M_+(\kappa)=\frac{q}{(x+\i)^\varsigma}H_+$\; and \;$\ran M_+(\kappa)=P_{H_+}\big(\frac{P}{Q_>Q_<(x+\i)^\varsigma}H_+\big)$
\item[\emph{(b)}] $\ker M_+(\kappa)=\Big\{\frac{Q_\leqslant}{P_<}r: \;r\textrm{ polynomial with }\deg r< \min\{\deg P_<-\deg Q_\le, \;\deg Q_>-\deg P_\ge\} \Big\}$
\item[\emph{(c)}] $\big(\ran M_+(\kappa)\big)^\perp=\Big\{\frac{\tilde{Q}_<}{\widetilde{P}_<}r:  \;r\textrm{ polynomial with }\deg r< \deg P_>-\deg Q_>\Big\}$
\item[\emph{(d)}] the following statements are equivalent:
\begin{itemize}
\item[\emph{(1)}] $M_+(\kappa)$ is a Fredholm Operator
\item[\emph{(2)}] $\ran M_+(\kappa)$  closed 
\item[\emph{(3)}] $\deg Q\le \deg P$ and $P$ without real zeros  
\item[\emph{(4)}] $0\not\in \overline{\kappa(\R)}$    
\end{itemize}
\end{itemize}
Here the zeros of the polynomial $q$ are the real zeros of $Q$,   $\varsigma:=\max\{\deg q, \deg P-\deg Q+\deg q\}$. Moreover $P=P_<P_\ge$, where the zeros of $P_<$($P_>$) and $P_\ge$ are exactly the zeros of $P$ in the lower(upper) half-plane and in the closed upper half-plane, respectively. $\tilde{P}$ denotes the polynomial whose coefficients are the complex conjugates of $P$. Analogous notations concern $Q$.
\end{The}\\
{\it Proof.}   For the closeness of $M_+(\kappa)$ write $\frac{P}{Q}$ in the form $\frac{P}{Q}=\frac{P_0}{Q_0} +\frac{p}{q}$ with  polynomials $P_0,Q_0,p,q$ such that $Q=Q_0q$,   $Q_0$ has no real zeros, $q$ has only real zeros,    $\deg P_0< \deg Q_0$, and 
$p$ and $q$ have no common zeros and satisfy $\varsigma=\max\{\deg p,\deg q\}$.\\
\hspace*{6mm}
Since $\kappa_0:=\frac{P_0}{Q_0}\big|_\R$ is bounded, $M_+(\kappa_0)$ is bounded. It follows 
$M_+(\kappa)= M_+(\kappa_0) +M_+(\frac{p}{q})$ and it remains to show that $M_+(\frac{p}{q})$ is closed. Let $h_n\in\dom M_+(\frac{p}{q})$ such that  $(h_n)$ converges to some $h\in H_+$ and $\big(M_+(\frac{p}{q})h_n\big)$ converges to some $k\in H_+$. By (\ref{AL1})(c), $\frac{p}{q}h_n\in H_+$. Hence one has $h_n\to h$ and $\frac{p}{q}h_n\to k$ in $L^2(\R)$. Since $M(\frac{p}{q})$ is closed, $h\in \dom M_+(\kappa)$ and $k=M_+(\frac{p}{q})h$ follows.\\
\hspace*{6mm}
(a)  $\dom M_+(\kappa)$ is dense by (\ref{DBWHO})(b). Arguing as above it remains to show $\dom M_+(\frac{p}{q})=\frac{q}{(x+\i)^\varsigma}H_+$. Let $h\in H_+$. Then, by (\ref{AL1})(c), $\frac{q}{(x+\i)^\varsigma}h\in H_+$  and $\frac{p}{q}\frac{q}{(x+\i)^\varsigma}h\in H_+$
implying $ \frac{q}{(x+\i)^\varsigma}H_+\subset \dom M_+(\frac{p}{q})$. For the converse inclusion argue $g\in \dom M_+(\frac{p}{q})$ $\Rightarrow$ $g\in H_+$, $\frac{p}{q}g\in L^2$ $\Rightarrow$ $h:=\frac{p}{q}g\in H_+$ by (\ref{AL1})(c). Hence $g=\frac{q}{p}h$, whence $\frac{1}{p}h\in H_+$ by (\ref{AL1})(d). Since 
$\varsigma=\max\{\deg p,\deg q\}$ one infers $k:=\frac{(x+\i)^\varsigma}{p}h\in L^2$, whence $k\in H_+$ applying (\ref{AL1})(c) to $\frac{1}{p}h\in H_+$. This shows $g=\frac{q}{(x+\i)^\varsigma}k\in \frac{q}{(x+\i)^\varsigma}H_+$. --- Now the claim about $\ran M_+(\kappa)$ is obvious.\\
\hspace*{6mm}
(b) Check the implications: $h_0\in\ker M_+(\kappa)$ $\Leftrightarrow$ $h_0\in \dom M_+(\kappa)$, $M_+(\kappa)h_0=0$  $\Leftrightarrow$ $h_0\in H_+$, $\kappa h_0\in L^2(\R)$, and $\kappa h_0\in H_-$. According to (\ref{AL1})(b) this means $h_0=R/P$, where $R$ is a polynomial with $\deg R<\min\{\deg P,\deg Q\}$ and $R=P_\ge Q_\le r$, whence the claim.\\
\hspace*{6mm}
(c) Using (a) one has $h_0\in \big(\ran M_+(\kappa)\big)^\perp$ $\Leftrightarrow$
$0=\langle h_0,\frac{P}{Q_>Q_<(x+\i)^\varsigma}h\rangle=\langle \frac{\tilde{P}}{(Q_>)^\sim \,(Q_<)^\sim \,(x-\i)^\varsigma}h_0,h\rangle$ $\forall$ $h\in H_+$ $\Leftrightarrow$ $\frac{\tilde{P}}{(Q_>)^\sim \,(Q_<)^\sim \,(x-\i)^\varsigma}h_0 \in H_-$. By (\ref{AL1})(b) this means $h_0=R/\tilde{P}$, where $R$ is a polynomial with $\deg R<\min\{\deg \tilde{P},\deg \tilde{Q}-\deg q +\varsigma\}=\deg P$ and $R=\tilde{P}_\ge (Q_>)^\sim \,r=\tilde{P}_\ge \tilde{Q}_< r$, whence the claim.\\
\hspace*{6mm}
(d)  By the forgoing results (1)$\Leftrightarrow$(2) holds. Moreover (3)$\Leftrightarrow$(4) is quite obvious. So we turn to 
(2)$\Leftrightarrow$(3). \\
\hspace*{6mm}
Show first  (3)$\Rightarrow$(2). Put $R:=\frac{P}{(Q_>)^\sim Q_< (x+\i)^\varsigma}$.  Since by the assumptions nominator and denominator of $R$ have equal degree and have no real zeros, $R$ and $\frac{1}{R}$ are bounded on $\R$. 
Hence $M(R)$ is a homeomorphism on $L^2(\R)$, whence $RH_+$ is closed. Since $(Q_>)^\sim=\tilde{Q}_<$ it follows by (\ref{AL1})(c) that $RH_+\subset H_+$. Hence $RH_+=\ran M_+(\kappa')$ for $\kappa':=\frac{P}{Q'}$ with $Q':=\tilde{Q}_< Q_< q$, whence $\dim (RH_+)^\perp<\infty$ by (c). So it suffices to show that $RH_+\subset \ran M_+(\kappa)$. 
By (\ref{AL1})(c), $\frac{Q_>}{\widetilde{Q}_< }H_+\subset H_+$. Hence, by (a), $\ran M_+(\kappa)\supset P_{H_+}\big(\frac{P}{Q_>Q_<(x+\i)^\varsigma} \frac{Q_>}{(Q_>)^\sim }H_+\big)=P_{H_+}(RH_+)=RH_+$.\\
\hspace*{6mm}
Now turn to (2)$\Rightarrow$(3). Consider first the case $\deg P<\deg Q$. Put $R:=Q_>Q_<(x+\i)^\varsigma$. Note $\deg R=\deg Q$ and $\ran M_+(\kappa)=P_{H_+} \big(\frac{P}{R}H_+\big)$ by (a). By (\ref{AL1})(c) and  since $M(\frac{R_>}{\widetilde{R}_< })$ is unitary on $L^2(\R)$, $\frac{R_>}{\widetilde{R}_< }H_+$ is a closed subspace of $H_+$. By (c) one has $\big(\frac{R_>}{\widetilde{R}_< }H_+\big)^\perp=\big\{\frac{r}{\widetilde{R}_<}:  \;\deg r< \deg R_>\big\}$ and hence 
$$\textrm{\SMALL{$\frac{P}{R}$}}H_+=\Big\{\textrm{\SMALL{$\frac{Pr}{R\widetilde{R}_<}$}}:  \;\deg r< \deg R_>\Big\}+\textrm{\SMALL{$\frac{P}{R_<\widetilde{R}_<}$}}H_+$$  
Applying (c) for $\kappa=\frac{P}{R_<\widetilde{R}_<}$ and $\kappa=
\frac{P_>}{\widetilde{P}_< }$ check $\big(\frac{P}{R_<\widetilde{R}_<}H_+\big)^\perp=
\big(\frac{P_>}{\widetilde{P}_<}H_+\big)^\perp=\big\{\frac{r}{\widetilde{P}_<}:  \;\deg r< \deg P_>\big\}$. Hence $\frac{P}{R_<\widetilde{R}_<}H_+$ is dense in $\frac{P_>}{\widetilde{P}_<}H_+$, which is closed in $H_+$. Note that $P_{H_+}\big\{\frac{Pr}{R\widetilde{R}_<}:  \;\deg r< \deg R_>\big\}\subset V:=\big\{\frac{B}{R_<\widetilde{R}_<}:\deg B<\deg R\big\}$.\\
\hspace*{6mm}
Now assume that $P_{H_+} \big(\frac{P}{R}H_+\big)$ is closed. Then the above considerations imply that $\frac{P_>}{\widetilde{P}_<}H_+ \subset V
+ \frac{P}{R_<\widetilde{R}_<}H_+$. Thus given $h_0\in H_+$ there is $h\in H_+$ and some polynomial $B$ with
$\deg B<\deg R$ such that $h=\frac{R_<\tilde{R}_<}{P_<\widetilde{P}_<}h_0+\frac{B}{P}$. Hence 
$|h|\ge \Big|\big|\frac{R_<\tilde{R}_<}{P_<\widetilde{P}_<}\big||h_0|-\big|\frac{B}{P}\big|\Big|$ with $m:=\deg R-\deg P_<-\deg P_> \ge 1$ and $l:=\deg B-\deg P\le m-1$. Choose $h_0\in H_+$ satisfying $|h_0|=(1+|x|)^{-3/4}$. ($h_0$ exists by (\ref{MHF}) since square-integrable and $\ln(1+|x|)\le \sqrt{2|x|}$.) 
Because of  $m-\frac{3}{4}-l>0$  the right side tends to $\infty$ like $|x|^{m-3/4}$ for $|x|\to \infty$ contradicting $h\in L^2(\R)$.\\
\hspace*{6mm}
To complete the proof of the implication (2)$\Rightarrow$(3) it remains to treat the case  that $\deg P\ge \deg Q$ and $P$ has a real zero. Proceeding as in the forgoing case, here one has  $\deg R=\deg P$. Hence, assuming that $P_{H_+} \big(\frac{P}{R}H_+\big)$ is closed,  one has $m\ge 1$ and  $l\le -1$ and the same contradiction follows.
\qed\\

From  (\ref{RSKCK}) one immediately obtains

\begin{Cor}\label{SWHFP} Let $\lambda\in\C$ and put $P^\lambda:=P+\lambda Q$. Then referring to $M_+(\kappa)$, $\lambda$ is 
\begin{itemize}
\item a Fredholm point  (i.e. $M_+(\kappa)-\lambda I$ is a Fredholm operator) iff $\lambda\not\in \overline{\kappa(\R)}$; if $\lambda$ is a Fredholm point, then $\dim\ker(M_+(\kappa)-\lambda I)=\max\{0,\deg Q_>-\deg P^\lambda_>\}$,
 $\dim\ran(M_+(\kappa)-\lambda I)^\perp=\max\{0,\deg P^\lambda_>-\deg Q_>\}$, and $\operatorname{ind}(M_+(\kappa)-\lambda I)=\deg Q_>-\deg P^\lambda_>$
\item a regular value (i.e. $M_+(\kappa)-\lambda I$ is continuously invertible) iff $\lambda\not\in \overline{\kappa(\R)}$ and $\deg Q_>\le \deg P^\lambda_>$
\item  in the resolvent set iff $\lambda\not\in \overline{\kappa(\R)}$ and $\deg Q_>= \deg P^\lambda_>$
\item a spectral value iff $\lambda\in \overline{\kappa(\R)}$ or  $\deg Q_>\ne \deg P^\lambda_>$
\item in the point spectrum iff \,$\deg Q_\le <\deg P^\lambda_<$ and $\deg P^\lambda_\ge <\deg Q_>$
\item in the continuous spectrum (i.e.  $M_+(\kappa)-\lambda I$ is injective with dense not closed range)  iff $\lambda\in \overline{\kappa(\R)}$, $\deg P^\lambda_> \le\deg Q_>$ and either 
$\deg P^\lambda_<\le \deg Q_\le$ or $\deg Q_>\le \deg P^\lambda_\ge$
\item in the residual spectrum (i.e.  $M_+(\kappa)-\lambda I$ is injective with not dense range) iff $\deg Q_> <\deg P^\lambda_>$
\end{itemize}
\end{Cor}
 The characterization of the Fredholm points of $M_+(\kappa)$  and the fact that at a Fredholm point either the kernel or the deficiency space is trivial,  are familiar from Krein's theory \cite{K62} for the  case of integrable kernel.  $\lambda\mapsto \deg P^\lambda_>$ is locally non-decreasing  due to the continuity of the roots of a polynomial on its coefficients \cite{CC89}. On  $\C\setminus \overline{\kappa(\R)}$ it is even locally constant, since there
 $P^\lambda_>=P^\lambda_\ge$. Hence, besides $\operatorname{ind}(M_+(\kappa)-\lambda I)$, also
$\dim\ker(M_+(\kappa)-\lambda I)$ and
 $\dim\ran(M_+(\kappa)-\lambda I)^\perp$ are constant on the components of $\C\setminus \overline{\kappa(\R)}$.

\begin{The}\label{DIWHO} Let $\kappa$ be a real rational function. Let $p$, $q$ be  real polynomials without common zeros such that  $\kappa=\frac{p}{q}$ and put $$ p-\i q=Q_+Q_-$$
where the zeros of the polynomials $Q_+$ and $Q_-$ are all in the upper and lower half-plane, respectively.
Then $M_+(\kappa)$ is closed symmetric with $\dom M_+(\kappa)=\frac{q'}{(x+\i)^\varsigma}H_+$, where the zeros of $q'$ are the real zeros of $q$ and $\varsigma:=\max\{\deg q',\deg p-\deg q +\deg q'\}$. The deficiency  indices are
  $$n_\pm:=\dim \ran\big(M_+(\kappa) \mp \i I\big)^\perp=\deg Q_\pm-\deg q_<$$
and the deficiency spaces  are
$$\ran\Big(M_+(\kappa)- \i I\Big)^\perp=\big\{\frac{q_<}{\widetilde{Q}_+}r: r\text{ polynomial with }\deg r<\deg Q_+-\deg q_<\}$$
$$\ran\Big(M_+(\kappa)+ \i I\Big)^\perp=\big\{\frac{q_<}{Q_-}r: r\text{ polynomial with }\deg r<\deg Q_--\deg q_<\}$$
\end{The}\\
{\it Proof.}  For $\dom M_+(\kappa)$ see (\ref{RSKCK})(a). Adopting the notation of (\ref{RSKCK}) one has $M_+(\kappa)- \i I
=M_+\big(\frac{P}{ Q}\big)$ with $P:=p-\i q$ and $Q=q$. Obviously $P$ and $Q$ have no common zeros. In view of (\ref{RSKCK})(c) note $\tilde{Q}_<=q_<$, $\tilde{P}_<
=\widetilde{Q}_+$, $\deg P_>=\deg Q_+$, $\deg{Q}_>=\deg{q}_>=\deg{q}_<$. Hence the formula for $\ran\big(M_+(\kappa)- \i I\big)^\perp$
 holds, which implies $n_+=\max\{0,\deg Q_+-\deg q_<\}$. It remains to show $\deg q_<\le\deg Q_+$. The assertions about the  deficiency space for $-\i$ follow similarly.\\
\hspace*{6mm}
Let $m_t$ denote the number of zeros in the upper half-plane of $p-\i tq$ for $t>0$. $p-\i tq$ has no real zeros. So by  the continuity of the roots of a polynomial, $m_t$ is locally constant and hence constant $=\deg Q_+$. Obviously $m_t$ is also the number of zeros in the upper half-plane of $\frac{1}{t}p-\i q$. Then by  the same continuity  the zeros of $q$ in the upper half-plane stay there for all $t$ large enough, whence $m_t\ge \deg q_>$. This yields the result.
 \qed\\


$W_\kappa$ for $\kappa$ in (\ref{DIWHO}) has unequal deficiency indices and hence no self-adjoint extension if $\max\{\deg p,\deg q\}$ is odd. 
 The deficiency indices $(n_+,n_-)$ of $W(x^l)$, $l\in\mathbb{Z}$, are  $(\frac{|l|}{2},\frac{|l|}{2})$ if $l$ is even and  $(\frac{|l+1|}{2},\frac{|l-1|}{2})$ otherwise.
  Other interesting examples are $W\big(\frac{x^2+1}{x}\big)$ and  
 $W\big(\frac{x^2-1}{x}\big)$ with deficiency indices $(1,1)$ and $(2,0)$, respectively. Compare the later with  $W\big(\frac{x^2-1}{x-2}\big)$ having deficiency indices $(1,1)$. --- This section is concluded by a much needed
 
\begin{Exa}\label{ENCWH}    There are  WH, even essentially self-adjoint semibounded ones, which are not closed, as for instance $W_{|x|}$ (or $W_{1/|x|}$).\\
\hspace*{6mm}
As to the proof, by (\ref{RSKCK})(a), $\dom M_+(|x|)=\frac{1}{x+\i}H_+$ and hence $R_\pm:=\ran (M_+(|x|)\mp \i I)= \frac{|x|\mp\, \i}{x+\i}H_+$. It suffices to show that $R_+$ is dense and $\ne H_+$.\\
\hspace*{6mm}
The latter is easily inferred. Assume the contrary. Then there is $h\in H_+$ satisfying $\frac{|x|- \i}{x+\i}h=\frac{1}{x+\i}$. Hence $h=\frac{1}{|x|-\i}\in H_+$ and $h=\check{h}\in H_-$, whence the contradiction $h=0$.\\
\hspace*{6mm} 
As to the former claim, let $h_\pm\in R^\perp_\pm$. Then $\langle h_\pm, \frac{|x|\mp\, \i}{x+\i}h\rangle=0$ for all $h\in H_+$, whence $k_\pm:=\frac{|x|\pm\, \i}{x-\i}h_\pm\in H_-$. Let $\phi_\pm$ be holomorphic on the upper half-plane with $\phi_{\pm,y} \to h_\pm$ for $0<y\to 0$ pointwise a.e. and in the mean (cf.\,(\ref{HS})). Similarly let $\psi_\pm$ be holomorphic in the lower half-plane with $\psi_{\pm,y}\to k_\pm$ for  $0>y\to 0$. Consider $\phi:=(z+\i)\phi_-\phi_+$ and $\psi:=(z-\i)\psi_-\psi_+$ holomorphic on the upper and lower half-plane, respectively. Check that $\phi_y$ for $0<y\to 0$ and $\psi_y$  for  $0>y\to 0$ converge pointwise a.e. and in $L^1_{loc}$ to $f:=(x+\i)h_-h_+$. By \cite[Theorem II]{C44} there is an entire function $\chi$ extending  $\phi$ and $\psi$. Let $\chi_\pm$ be equal to $\phi_\pm$ and $\psi_\pm$ on the upper and lower half-plane, respectively. 
Let $z=x+\i y$, $|z|>1$. Then $\pi |\chi(z)|\le\int_{|w|\le 1}|\chi(z+w)|\d^2w\le (|z|+2)\big(\int_{|w|\le 1}|\chi_+(z+w)|^2\d^2w\,\int_{|w|\le 1}|\chi_-(z+w)|^2\d^2w\big)^{1/2}$. Now $\int_{|w|\le 1}|\chi_\pm(z+w)|^2\d^2w\le \int_{-1}^1\int_{-\infty}^\infty|\chi_\pm\big(u+\i(y+v)\big)|^2\d u\d v\le\int_{-1}^1C\d v =2C $ for some finite constant $C$. Hence $|\chi(z)|\le C'  |z|$ $\forall \,|z|>1$. So $\chi$ is a polynomial $a+bz$. For $y>0$, $x\mapsto \phi_{-,y}(x)\phi_{+,y}(x)=\frac{a+bz}{z+\i}$ is integrable on $\R$. This implies $a=b=0$ and hence $\chi=0$. Therefore either $h_+=0$ or $h_-=0$. This means that one and hence both deficiency spaces are $\{0\}$. So $R_+$ is dense.
\end{Exa}

\section{Semibounded Wiener-Hopf operators}\label{SBWHO}

In (\ref{RSBWHO}),\,(\ref{SBS})  a semibounded densely defined WH $W_\kappa$ is  expressed in a canonical way by the product of a closable operator and its adjoint. Replacing the operator by its closure one obtains a self-adjoint extension $\tilde{W}_\kappa$ of  $W_\kappa$, which is semibounded by the same bound. The bound is not an eigenvalue of the extension. $\tilde{W}_\kappa$   is shown to be the Friedrichs extension of 
$W_\kappa$.

\begin{Lem} \label{CPEFIPS}  Let $\gamma:\R\to\C$ be measurable. Put $E:=\gamma^{-1}(\C\setminus\{0\})$ and let  $A:=P_+\mathcal{F}M(\gamma)P_E^*$.  Then $A$ is densely defined and  $A^*=P_EM(\overline{\gamma})\F^{-1}P_+^*$ holds, and  $\dom A^*=\dom W_{\overline{\gamma}}$ is either  $\{0\}$ or  dense. If $\gamma$ is not almost zero, then
$A^*$ is injective. If $\gamma\ne 0$ not a.e., then $A$ is injective. If $\gamma\ne 0$ a.e., then $\ker A=
\{f\in L^2(\R): \gamma f\in H_-\}$, 
which equals $\{0\}$ if and only if $\frac{1}{\gamma}$ is not proper. Finally  $\ker A^*A=\ker A$, $W(|\gamma|^2)=AA^*$, and
 \begin{equation*}
\gamma  \textrm{ proper } \Leftrightarrow \dom W(|\gamma|^2) \ne \{0\} \Leftrightarrow \dom A^*\ne\{0\} \Leftrightarrow A \textrm{ closable}\tag{$\star$}
 \end{equation*}
\end{Lem}\\
{\it Proof.} Note  $M(\gamma)P_E^*=P_E^*M(\gamma|_E)$ and $P_EM(\gamma)=M(\gamma|_E)P_E$, where $M(\gamma|_E)$ denotes the multiplication operator in $L^2(E)$. Hence $A=
P_+\mathcal{F} P_E^*M(\gamma|_E)$ is densely defined, and $A^*=P_EM(\overline{\gamma})\F^{-1}P_+^*$  by \cite[13.2(2)]{R74}. Note that $\dom A^*=\dom M(\overline{\gamma})\F^{-1}P_+^*=\dom W_{\overline{\gamma}}$, whence the claim on  $\dom A^*$ by (\ref{DBWHO})(i). \\
 \hspace*{6mm} 
 First let  $E$ be proper. Then $A^*g
=M(\overline{\gamma}|_E)\big(P_E\mathcal{F}^{-1}P_+^*g\big)=0$ implies $P_E\mathcal{F}^{-1}P_+^*g=0$, whence $g=0$. Similarly, $Ak=
P_+\mathcal{F} \big( P_E^*M(\gamma|_E)k\big)=0$ means that $f_-:=P_E^*M(\gamma|_E)k\in H_-$ and vanishes on $\R\setminus E$, whence $f_-=0$ and hence $k=0$.   ---  Now assume at once  $E=\R$. Then $A^*=M(\overline{\gamma})\F^{-1}P_+^*$ is injective as $M(\overline{\gamma})$ is injective. Furthermore, $f\in \ker A\Leftrightarrow f, \gamma f\in L^2(\R)$ with $P_+\F(\gamma f)=0 \Leftrightarrow f\in L^2(\R),  \gamma f\in H_- \Leftrightarrow  
  f\in L^2(\R),  \overline{f}\in\frac{1}{\overline{\gamma}} H_+$,   
 and recall that $\frac{1}{\gamma}$ is proper $\Leftrightarrow$
 $\frac{1}{\overline{\gamma}} h\in L^2(\R)$  for some $h\in H_+\setminus\{0\}$.\\
\hspace*{6mm}
 $\ker A^*A=\ker A$ is obvious since either  $\gamma=0$ a.e. or $A^*$ is injective.  $M(\gamma)P_E^*P_EM(\overline{\gamma})=M(\gamma)M(1_E)M(\overline{\gamma})=M(|\gamma|^2)$, whence  $AA^*=W(|\gamma|^2)$.\\
 \hspace*{6mm}
Turn to the final claim ($\star$). Recall that  $\dom A^*=\dom W_{\overline{\gamma}}$ is either trivial or dense. So the last equivalence is standard and the remaining equivalences  hold by (\ref{DBWHO})(f). \qed\\
  
 In view of  (\ref{CPEFIPS})($\star$) recall  the results on the domain of a WH in (\ref{DBWHO}).    If $A$ is closable it   need not be closed, even if $AA^*$ is closed. $\big($Indeed,   $W_{x^2}=AA^*$  is closed and $\ne \tilde{W}_{x^2}=\overline{A}A^*$ by (\ref{DIWHO}),(\ref{RSBWHO}).$\big)$ Recall that $\ran A^*$ is not dense if $A$ is not injective.
$A$ is not injective for $\gamma\ne0$ a.e. if for instance $\frac{q}{\gamma}\F^{-1}s\in L^2(\R)$ with $q$ a polynomial and $s$ a Schwartz function  with support in $]-\infty,0]$. (Indeed,  $q\F^{-1}s$ is a Schwartz function in $H_-$, whence $\frac{q}{\gamma}\F^{-1}s  \in \ker A$.) Recall that $A$ is injective  if and only if  $\gamma\ne 0$ not a.e. or  $\gamma\ne 0$ a.e. and $\frac{1}{\gamma}$ not proper.

 \begin{Lem}\label{AL2} Let $\gamma$ be proper. Suppose that
   $\gamma\ne 0$ not a.e. or that $\gamma\ne 0$ a.e. and $\frac{1}{\gamma}$ is not proper.  Then $A$ is closable and $\overline{A}$ is injective.
 \end{Lem}\\
{\it Proof.} $A$ is closable  by (\ref{CPEFIPS})($\star$). According to (\ref{DBWHO})(e)($p=\infty$) one has 
$|\gamma| \e^j\le 1$ for some $j$. Put $j':=1_{\{j\le 0\}}j-|x|^{1/2}$. Then $|\gamma| \e^{j'}\le 1$, and $\e^{j'}$ and $|\gamma|\e^{j'}$ are square integrable. By (\ref{MHF}) there is $h\in H_+$ with $|h|=\e^{j'}$. Hence $h\in \dom M_+(\overline{\gamma})\setminus\{0\}$. 
Then $hD\subset \dom M_+(\overline{\gamma})$ for $D:=\frac{1}{x+\i} H_+^\infty$. $D$ is dense in $H_+$ since
$D=\Gamma(H^\infty(\mathbb{T}))$ (see (\ref{RTO})). --- 
Now let $f\in\ker\overline{A}=
(\ran A^*)^\perp$. For $d\in D$ one has $0=\langle f,P_EM(\overline{\gamma})P_{H_+}^*(hd)\rangle =\langle P_E^*f,\overline{\gamma}hd\rangle=\langle \gamma\overline{h}P_E^*f,d\rangle$, whence $\gamma\overline{h}P_E^*f\in H_-$. 
If $\gamma\ne0$ not a.e., then it follows $\gamma\overline{h}P_E^*f=0$, whence $f=0$.  If $\gamma\ne 0$ a.e., then $h\overline{f}=\frac{1}{\overline{\gamma}}h'$ for some $h'\in H_+$. Assume $f\ne 0$. Then $h'\ne 0$ and $|\gamma|^{-1}||h'|$ is integrable. Hence  $|\gamma|^{-1}$ is proper by (\ref{MHF}) and  (\ref{DBWHO})(e)($p=1$), which however  is excluded by the premise.\qed

The foregoing lemma is  needed only in sec.\,\ref{ISIO}. The main result of this section follows.

\begin{The}\label{RSBWHO}  Let $\kappa\ge 0$. Put $E:=\kappa^{-1}(\R\setminus\{0\})$ and $A:=P_+\F M(\sqrt{\kappa})P_E^*$. Then  
$W_\kappa = AA^*$.\\
\hspace*{6mm}
 Now let $\kappa$ be proper not almost zero. Put $ \tilde{W}_\kappa:= \overline{A}A^*$. Then 
\begin{itemize}
\item[\emph{(a)}] $W_\kappa$  is  densely defined symmetric nonnegative and $\tilde{W}_\kappa$ 
  is an injective nonnegative self-adjoint extension of $W_\kappa$.  
\item[\emph{(b)}] 
$\dom W_\kappa$ is a core of $A^*$ and $\dom A^*\,\cap\, \ran(I+W_\kappa)^\perp=\emptyset$ holds.
\item[\emph{(c)}] 
 $ \tilde{W}_\kappa$ is the  Friedrichs extension of $W_\kappa$. 
\end{itemize}  
\end{The}
{\it Proof.} (a) Apply (\ref{CPEFIPS}) for $\gamma:=\sqrt{\kappa}$. Accordingly,  $W_\kappa=AA^* $ and, if 
 $\kappa$ is proper,  $W_\kappa$ is densely defined and symmetric nonnegative by (\ref{DBWHO})(n),(o), and
  $A^*$ is densely defined.  Then $\overline{A}=A^{**}$ and
by  \cite[13.13(a)]{R74} $\overline{A}A^*$ is self-adjoint. Clearly $\overline{A}A^*$ is nonnegative. Check that $\overline{A}A^*$ is injective  as  $A^*$ is injective  
by  (\ref{CPEFIPS}) for $\kappa$ not almost zero.\\
\hspace*{6mm}
 (b),\,(c) According to \cite[Theorem X.23]{RS75}, 
$ \tilde{W}_\kappa$  is the  Friedrichs extension only if 
$\dom  \tilde{W}_\kappa \subset H_{W_\kappa}$, where $H_{W_\kappa}$ is the completion of $\dom W_\kappa$ with respect to the sesquilinear form $\langle g,g'\rangle_{W_\kappa}:=\langle g,g'\rangle+\langle g,W_\kappa g'\rangle$.\\
\hspace*{6mm}
Endow $\dom A^*$ with the inner product $\langle g,g'\rangle_{A^*}:=\langle g,g'\rangle+\langle A^*g,A^*g'\rangle$, by which $\dom A^*$ becomes a Hilbert space $\mathcal{K}$ since $A^*$ is closed. Then the  subspace $\dom \overline{A}A^*$ is dense in $\mathcal{K}$
since $\dom \overline{A}A^*$ is a core for $A^*$, see  \cite[13.13(b)]{R74}. One easily checks that  $H_{W_\kappa}$ 
 is the closure of $\dom W_\kappa$ in $\mathcal{K}$. Therefore $\dom \overline{A}A^*\subset H_{W_\kappa}$ if and only if $H_{W_\kappa}=\mathcal{K}$, which means that  $\dom W_\kappa$ is a core of $A^*$. 
 A short computation shows also that $H_{W_\kappa}=\mathcal{K}$ is equivalent to
 $\dom A^*\,\cap\, \ran(I+W_\kappa)^\perp=\emptyset$.\\
\hspace*{6mm}
Hence it remains to show  $\dom A^*\,\cap\, \ran(I+W_\kappa)^\perp=\emptyset$. Explicitly this means  that $h_0=0$ if 
\begin{equation*}
h_0\in H_+ \textrm{ satisfies } \sqrt{\kappa}h_0\in L^2(\R) \textrm{ and } \langle h_0,(1+\kappa)h\rangle =0 \;\forall\; h\in H_+ \textrm{ with } \kappa h\in L^2(\R) \tag{$\star$}
\end{equation*} 
To this end a sequence $(\alpha_n)_n$ in $H^\infty_+$ will be constructed with
$|\alpha_n|=\e^{j_n}$ a.e.,  where $j_n(x):=0$ if $1+\kappa(x)< n$ and $j_n(x):=-\frac{1}{2}\ln\big(1+\kappa(x)\big)$ otherwise, and satisfying $\alpha_{n_k}\to 1$ pointwise a.e. for some subsequence $(n_k)$. Provided $(\alpha_n)_n$  set $h_n:=\alpha_nh_0$. Then  $h_n\in H_+$ and almost everywhere 
$\big(1+\kappa(x)\big)|h_n(x)|$ is less than $n|h_0(x)|$ if $1+\kappa(x)\le n$ and equals $\sqrt{1+\kappa(x)}|h_0(x)|$ otherwise, which proves $(1+\kappa)h_n\in L^2(\R)$. Moreover $\sqrt{1+\kappa}|h_n|\le\sqrt{1+\kappa}|h_0|$ since $|\alpha_n|\le 1$ and $\sqrt{1+\kappa}h_{n_k}\to \sqrt{1+\kappa}h_0$ pointwise a.e., whence $\sqrt{1+\kappa}h_{n_k}\to \sqrt{1+\kappa}h_0$ in $L^2(\R)$ by dominated convergence. Thus ($\star$) holds for $h=h_{n_k}$, whence $0= \langle h_0,(1+\kappa)h_{n_k}\rangle=\langle \sqrt{1+\kappa}h_0,\sqrt{1+\kappa}h_{n_k}\rangle\to \langle \sqrt{1+\kappa}h_0,\sqrt{1+\kappa}h_0\rangle=\norm{\sqrt{1+\kappa}h_0}^2$ implying $h_0=0$.\\
\hspace*{6mm}
We turn to the construction of $(\alpha_n)_n$. By (\ref{DBWHO})(g), 
$\frac{\ln(1+\kappa)}{1+x^2}$  and hence all $\frac{j_n}{1+x^2}$ are integrable. For convenience we pass from $\R$ to the torus $\mathbb{T}$ by means of the Cayley transformation $C$ (see (\ref{RTO})). So let $\tilde{j}_n:=j_n\circ C^{-1}$, which is integrable on $\mathbb{T}$. Put
$$ F_n(w) :=\frac{1}{2\pi}\int_0^{2\pi}\frac{\e^{\i t}+w}{\e^{\i t}-w}\,\tilde{j}_n(\e^{\i t})\d t$$
for $w\in \mathbb{D}$. Then $\exp\circ F_n$ is an  outer function. Let  $\tilde{\alpha}_n$  denote its nontangential limit a.e.\;on $\mathbb{T}$. It satisfies $|\tilde{\alpha}_n|=\e^{\tilde{j}_n}$ a.e. Hence  
$\tilde{\alpha}_n\in H^\infty(\mathbb{T})$. It remains to show the existence of a subsequence $(n_k)$ satisfying $\tilde{\alpha}_{n_k}\to 1$ a.e. The formula $Tf(z):=\lim_{r\uparrow 1}\frac{1}{2\pi}\int_0^{2\pi}\frac{\e^{\i t}+rz}{\e^{\i t}-rz}\,f(\e^{\i t})\d t$, $z\in\mathbb{T}$ defines a bounded operator on $L^1(\mathbb{T})$ into weak-$L^1(\mathbb{T})$, whence $|\{z\in \mathbb{T}: |Tf|>\delta\}|\le C\norm{f}_1/\delta$
for all $\delta>0$ and $f\in L^1(\mathbb{T})$. Therefore, if $f_n\to 0$ in $L^1(\mathbb{T})$, then $Tf_n\to 0$ in probability,
which implies $Tf_{n_k}\to 0$ a.e. for some subsequence $(n_k)$. This applies to $(\tilde{j}_n)_n$ yielding $\tilde{\alpha}_{n_k}=\exp\circ T\tilde{j}_{n_k}\to 1$ a.e.\qed\\

For  $\kappa\ge 0$ and $W_\kappa$  densely defined recall that the deficiency subspace $\ran(I+W_\kappa)^\perp$ of $W_\kappa$ at $-1$ is  trivial if and only if $W_\kappa$ is  essentially self-adjoint.

\begin{SBS}\label{SBS} Let $\kappa$ be real semibounded. Then there are $\alpha>0$ and $\eta\in\{1,-1\}$ such that $\varkappa:=\alpha 1+\eta \kappa \ge 0$. Clearly $W_\varkappa=\alpha I +\eta W_\kappa$.  Let $W_\kappa$ be densely defined. Then so is $W_\varkappa$, and  according to (\ref{RSBWHO}) there is the injective nonnegative self-adjoint extension $\tilde{W}_\varkappa$ of $W_\varkappa$. So 
\begin{equation}\label{SBSE}
\tilde{W}_\kappa:= -\eta\alpha I+\eta \tilde{W}_\varkappa
\end{equation}
is a semibounded self-adjoint extension of $W_\kappa$ with bound $-\eta \alpha$, which is not an eigenvalue of $\tilde{W}_\kappa$. It is the Friedrichs extension of $W_\kappa$.
\end{SBS}

\section{Isomorphic Singular Integral Operators}\label{ISIO}

This section is concerned with the  symmetric singular integral operator  in $L^2(E)$ for proper $E$ or $E=\R$ 

\begin{equation}\label{OSIO}
\big(L_\phi f\big)(x):=\textrm{\SMALL{$\frac{1}{2}$}}\phi(x)f(x)+\textrm{\SMALL{$\frac{1}{2\pi\operatorname{i}}$}}\int_E\frac{\sqrt{\phi(x)}\sqrt{\phi(y)}}{y-x}f(y)\operatorname{d}y
\end{equation}
(in the sense of the principal value at $x$) where  $\phi:E\to\R$ is measurable positive. $L_\phi$ will turn out to be closely related to   $W_\kappa$, where $\kappa$ extends $\phi$ on $\R$ by zero.\\
\hspace*{6mm}
$L_\phi$ belongs to the  studied class of singular integral operators  in $L^2(E)$ of type Hilbert transformation 
\begin{equation}\label{GHTLSIO} 
(L(a,b)f\big)(x):=a(x)f(x)+\frac{1}{\i \pi}\int_E b(x)\overline{b(y)} (y-x)^{-1}f(y)\d y
\end{equation}
 where $a,b$ are measurable functions  on $E$ with $a$ real and $b\ne 0$ a.e. There is the obvious   unitary equivalence  
\begin{equation}\label{GCSIO}
U L(a,b)\, U^{-1}\; =\; L_\phi +M(\alpha) 
\end{equation}
for $\phi=2|b|^2$ and $\alpha=a-|b|^2$, where $U$ is the multiplication operator by $\overline{b}/|b|$ and $M(\alpha)$ the multiplication operator by $\alpha$ in $L^2(E)$. So we are concerned with the case $a=|b|^2$.\\
\hspace*{6mm}
The operator $L(a,b)$  for bounded $b$ and bounded below $a$ is treated
 by  Rosenblum in \cite{R66}. It is shown to be
self-adjoint on $\dom L(a,b)=\dom M(a)$,  and is diagonalization is achieved. See also \cite{Pu69}  and the literature cited in  \cite{R66},   \cite{Pu69}.  
The really  unbounded case however is there when $b$ is unbounded.
 \cite{P73}  is concerned with this case replacing $L(a,b)$ by the limit of truncated $L(a_n,b_n)$ which are bounded. Our analysis of $L_\phi$ will show  (\ref{TSIO}) that $L(a,b)\subset M(\alpha)$ if the extension of $b$ on $\R$ by zero is not proper. 
 Hence for $L(a,b)$ in (\ref{GHTLSIO}) being not trivial it is necessary that the extension of $b$ is proper. In this case $L_\phi$ in (\ref{GCSIO}) has a self-adjoint extension (\ref{CKPAE}).
 
The Hilbert transformation $H$ on $L^2(\R)$ is defined by the singular integral 

\begin{equation}
Hf(x)= \textrm{\SMALL{$\frac{1}{\operatorname{i}\pi}$}}\int_{-\infty}^\infty\frac{f(y)}{y-x}\operatorname{d}y
\end{equation}
Recall its representation 
\begin{equation}\label{SRHT}
H=\F^{-1}M(\sgn)\F=-\F M(\sgn)\F^{-1}
\end{equation} 
on $L^2(\R)$ with $\sgn$ the signum function on $\R$ (see e.g.\,\cite[Theorems 91,\,95]{T48} or \cite[Lemma 1.35]{D79} or \cite[Chapter II 4.3]{R05}, and for more details \cite[Teorema 1.1.1]{T73}). Let us introduce 
\begin{equation}\label{THTE}
H_E:=P_E\,H\,P_E^*
\end{equation} 
 the trace on $L^2(E)$ of the Hilbert  transformation $H$. Its spectrum is determined in  \cite{W60}. For  $E=[a,b]$, $-\infty\le a<b\le\infty$, $H_E$ is called finite and semi-finite Hilbert transformation if $E$ is bounded and semi-bounded, respectively. Its spectral representation is achieved in \cite{KP59}.

 \subsection{Isometry relating $\tilde{L}_{\phi}$ to $\tilde{W}_\kappa$}

In what follows we use the polar decomposition $C=S|C|$ of a closed densely defined operator $C$ from a Hilbert space $\mathcal{H}$ into another $\mathcal{H}'$ (see e.g. \cite[VIII.9]{RS75}).  $|C|$ denotes the square root of the self-adjoint nonnegative operator $C^*C$ in $\mathcal{H}$. One has $\dom |C|=\dom C$ and $\dom C^*C$ is a core for $C$.  $S$ is a partial isometry from $\mathcal{H}$ into  $\mathcal{H}'$.  Its initial space $(\operatorname{ker}S)^\perp$ equals $\overline{\ran C^*}=\overline{\ran |C|}=\overline{\ran |C|^2}$. Similarly, its final space $\ran S$ equals
$\overline{\ran C}=\overline{\ran |C^*|}=\overline{\ran |C^*|^2}$. The partial isometry $S^*$ satisfies $(\ker S^*)^\perp = \ran S$ and $\ran S^*=(\ker S)^\perp$. Important is the relation 
\begin{equation}\label{FRS}
CC^*=SC^*CS^*
\end{equation}
which means that  the reductions of $CC^*$ and $C^*C$ on the orthogonal complements of their  respective null spaces are Hilbert space isomorphic by the restriction of the partial isometry $S$ to its initial state $(\ker S)^\perp$ and  its final state $\ran S$.\\

Throughout this section  $\kappa\ge 0$ and $\phi>0$ are related to each other by $\phi=\kappa|_E$ for $E=\kappa^{-1}(\R_+)$.

\begin{Lem}\label{MLEM}  Let $\kappa\ge 0$. Recall $A=P_+\mathcal{F}M(\sqrt{\kappa})P_E^*$, $A^*=P_E M(\sqrt{\kappa})\F^{-1}P_+^*$. Then 
\begin{itemize}
\item[\emph{(a)}] $W_\kappa=AA^*$
\item[\emph{(b)}] $L_{\phi}=\frac{1}{2}M(\phi) +\frac{1}{2}
M(\sqrt{\phi})H_EM(\sqrt{\phi})=P_EM(\sqrt{\kappa})P^*_{H_+}P_{H_+}M(\sqrt{\kappa})P_E^* = A^*A$
\end{itemize}
Now suppose that $\kappa$ is proper. 
Let $\overline{A}=T|\overline{A}|$ be the polar decomposition of $\overline{A}$.
Then the  partial isometry  $T:L^2(E)\to L^2(\R_+)$ is surjective, its adjoint $T^*$ is injective, and
\begin{itemize}
\item[\emph{(c)}] $\tilde{W}_\kappa:=\overline{A}A^*$ is a self-adjoint extension of $W_\kappa$. 
\item[\emph{(d)}]
$\tilde{L}_{\phi}:=A^*\overline{A}$ is a self-adjoint extension of $L_{\phi}$.
\item[\emph{(e)}] $\ker \tilde{L}_{\phi}=\ker\overline{A}$, $\big(\ker T)^\perp=\big(\ker \tilde{L}_{\phi}\big)^\perp=\overline{\ran A^*}$.
\end{itemize}
\end{Lem}
{\it Proof.} (a) See (\ref{RSBWHO}). (b)  $A^*A=P_EM(\sqrt{\kappa})\mathcal{F}^{-1}P_+^*P_+\mathcal{F}M(\sqrt{\kappa})P_E^*=P_EM(\sqrt{\kappa})P^*_{H_+}P_{H_+}M(\sqrt{\kappa})P_E^* =P_EM(\sqrt{\kappa})\frac{1}{2}(I+H)M(\sqrt{\kappa})P_E^*$ by (\ref{SRHT}). Hence $A^*A=L_{\phi}$ as $P_EM(\sqrt{\kappa})=M(\sqrt{\phi}) P_E$, $M(\sqrt{\kappa})P_E^*=P_E^*M(\sqrt{\phi})$. Recall (\ref{THTE}). --- (c) and (d) are obvious. --- Finally, $T$ is surjective as $\operatorname{ran}T=(\,\operatorname{ran}\overline{A}\,)^-=(\operatorname{ker}A^*)^\perp=L^2(\R_+)$ by (\ref{CPEFIPS}), and (e) follows directly  from  polar decomposition. 
\qed\\

\hspace*{6mm} 
As a first result we note that $L_\phi$ in (\ref{OSIO}) can be trivial. 

\begin{Cor}\label{NSOL}  Let $\kappa\ge 0$ be not proper. Then $L_\phi\subset 0$. Moreover, $\dom L_\phi=\{0\}$ if and only if either $\kappa\ne 0$ not a.e.  or  otherwise $\frac{1}{\kappa}$ is not proper. If $\kappa\ne 0$ a.e. and $\frac{1}{\kappa}$ is proper, then $\dom L_\phi= \frac{1}{\sqrt{\phi}}H_-\cap L^2(\R)$. 
 \end{Cor}\\
{\it Proof.}  $\dom A^*=\{0\}$ by  (\ref{CPEFIPS})($\star$), whence $L_\phi\subset 0$ and $\dom L_\phi=\ker A$ because of  (\ref{MLEM})(b), (\ref{CPEFIPS}). Recall (\ref{DBWHO})(f) by which $1/\kappa$ is not proper  if and only if $1/\sqrt{\kappa}$ is not proper. The remainder  is shown in  (\ref{CPEFIPS}).\qed
  
For example $\kappa=\e^x$ is positive and,  by (\ref{DBWHO})(h),  $\kappa$  and $1/\kappa$ are not proper. Hence $\dom L_{\e^x}=\{0\}$. \\
\hspace*{6mm}
 Now (\ref{NSOL}) and (\ref{GCSIO}) imply that for $L(a,b)$   in (\ref{GHTLSIO}) being not trivial it is necessary that the extension of  $b$ on $\R$ by zero be proper:
 
\begin{The}\label{TSIO} Let the extension of  $b$ on $\R$ by zero  be not proper. Then $L(a,b)\subset M(a-|b|^2)$ with
\begin{itemize}
\item[\emph{(a)}] $\dom L(a,b)=\big\{f\in L^2(\R): |b|f\in H_-,   (a-|b|^2)f\in L^2(\R)\big\}$ if $E=\R$
\item[\emph{(b)}] $\dom L(a,b)=\{0\}$ if and only if $E$ is proper, or  $E=\R$  and $\frac{1}{b}$ is not proper,  or  $E=\R$ and $\frac{1}{b}$ is  proper and $\frac{a-|b|^2}{b}$ is not proper
\end{itemize}
\end{The}

\begin{Cor}\label{SBARZ} 
Let $E=\R$ and suppose that $\frac{1}{\phi}$ is bounded.. Then $L_\phi$ is densely defined and
\begin{equation*}
\dom L_\phi=\{\textrm{\SMALL{$\frac{1}{\sqrt{\phi}}$}}(h+k):h\in H_+, \sqrt{\phi}h\in L^2(\R), \,k\in H_-\},  \quad L_\phi \big(\textrm{\SMALL{$\frac{h+k}{\sqrt{\phi}}$}}\big)=\sqrt{\phi}h\tag{$\star$}
\end{equation*}
with $\ker L_\phi=\frac{1}{\sqrt{\phi}}H_-$ and  $(\ker L_\phi)^\perp=\{\sqrt{\phi}h\in L^2(\R):h\in H_+\}$.\\
\hspace*{6mm}
 If $\phi$ is proper, then 
 $\ker \tilde{L}_\phi=\overline{\ker L_\phi}\ne L^2(\R)$. If $\phi$ is not proper, then $\overline{L}_\phi=0$.
\end{Cor} \\
{\it Proof.} Note $\kappa=\phi$, $L_\phi=M(\sqrt{\phi})P^*_{H_+}P_{H_+}M(\sqrt{\phi})$ by (\ref{MLEM})(b).
Then ($\star$) is easily verified and $\ker L_\phi=\frac{1}{\sqrt{\phi}}H_-$ follows. Check $f\in L^2(\R)$, $f\perp \frac{1}{\sqrt{\phi}}H_-$ $\Leftrightarrow$ $f\in L^2(\R)$, $0=
\langle f,\frac{1}{\sqrt{\phi}}k\rangle=\langle \frac{1}{\sqrt{\phi}}f,k\rangle$ $\forall$ $k\in H_-$ $\Leftrightarrow$  $f\in L^2(\R)$, $\frac{1}{\sqrt{\phi}}f\in H_+$, proving the claim on $(\ker L_\phi)^\perp$.\\
\hspace*{6mm}
Note $D:=\dom A^*=\{h\in H_+:\sqrt{\phi}h\in L^2(\R)\}$, whence  $\ran A^*=\{\sqrt{\phi}h\in L^2(\R):h\in H_+\}$ which equals $(\ker L_\phi)^\perp$. \\
\hspace*{6mm}
Now let $\phi$ be proper.  Then $D$ is dense  by (\ref{CPEFIPS}) and $\ran A^*\ne\{0\}$.
Hence by (\ref{MLEM})(e) on has $\ker \tilde{L}_\phi=(\ran A^*)^\perp=\overline{\ker L_\phi}\ne L^2(\R)$. Moreover, $\dom L_\phi$ is dense, since $f\in (\dom L_\phi)^\perp$ $\Rightarrow$ $0=\langle f,(h+k)/\sqrt{\phi}\rangle=\langle f/\sqrt{\phi},(h+k)\rangle$ $\forall$ $h\in D$, $k\in H_-$ $\Rightarrow$ $f/\sqrt{\phi}\in H_+$, $\langle f/\sqrt{\phi},h\rangle$ $\forall$ $h\in D$ $\Rightarrow$ $f/\sqrt{\phi}=0$ $\Rightarrow$ $f=0$. --- If $\phi$ is not proper then $D=\{0\}$ by (\ref{CPEFIPS}), whence  $\dom L_\phi=\frac{1}{\sqrt{\phi}}H_-$. The latter equals $\ker L_\phi$. Finally,   
$(\ker L_\phi)^\perp=\{\sqrt{\phi}h\in L^2(\R):h\in H_+\}=\{0\}$ since $D=\{0\}$.\qed

For example, $L_{\e^{|x|}}\subset 0$ with dense $\dom L_{\e^{|x|}}=\e^{-|x|/2}H_-$.\\

\hspace*{6mm} 
 The main outcome of this section is

 \begin{The}\label{CKPAE}  Let $\kappa \ge 0$ be proper. 
 Then
 \begin{equation}\label{FRWHSIO}
 \tilde{L}_\phi=\,T^*\;\tilde{W}_\kappa\,T
 \end{equation}
  Let $\tilde{L}_\phi'$ and $T'$ denote the reduction of $\tilde{L}_\phi$ and the restriction of $T$ 
 on the orthogonal complement of the null space of $\tilde{L}_\phi$, respectively.   Then $T'$ 
 is a Hilbert space isomorphism onto 
 $L^2(\R_+)$ satisfying
 $\tilde{L}_\phi'=\,T'^{-1}\;\tilde{W}_\kappa\,T'$.
   \end{The}\\
{\it Proof.}  (\ref{FRWHSIO})  follows immediately from polar decomposition $\overline{A}=T|\overline{A}|$ (cf.\,(\ref{FRS})). By (\ref{MLEM}) $T$ is surjective. Hence $T'$ is a Hilbert space isomorphism by  (\ref{MLEM})(e). Finally recall  (\ref{FRWHSIO}).\qed\\

Of course one has    $\tilde{W}_\kappa =T\tilde{L}_\phi T^*$ as well. It allows to study $\tilde{W}_\kappa$ starting from  $\tilde{L}_\phi$. For the still remarkable bounded case see (\ref{RBC}). By (\ref{SBS}) it is easy to extend (\ref{CKPAE}) to semibounded symbols.

\begin{Cor}\label{ACSIO} 
Let  $\frac{\kappa}{1+x^2}$ be integrable. Then $\tilde{L}_\phi'$ is absolutely continuous.
\end{Cor}\\
{\it Proof.} Recall that $\tilde{T}_{\kappa\circ \gamma}$ and hence $\tilde{W}_\kappa$ (see (\ref{RTO})) is absolutely continuous by \cite{R65}. Apply (\ref{CKPAE}).\qed\\

The special case of (\ref{ACSIO}) that $\phi$ is bounded  is treated in \cite[sec.\,3]{R66} and \cite[Theorem]{Pu69}. --- We remind that for $\tilde{L}_\phi$ being injective it is necessary that $E$ is proper or that $1/\kappa$ is not proper.

\begin{Cor} \label{CCKPAE} Let $\kappa\ge 0$ be proper. Suppose that  $E$ is proper or $1/\kappa$ is not proper. Then $\tilde{L}_\phi$ is injective and $T:L^2(E)\to L^2(\R_+)$ is a Hilbert space isomorphism with  $$\tilde{L}_\phi=\,T^{-1}\tilde{W}_\kappa\,T$$ 
 \end{Cor}
{\it Proof.}   $\overline{A}$ is injective by (\ref{AL2}). Hence $\{0\}=\ker \tilde{L}_\phi =\ker T$ and $T$ is an isomorphism. Apply (\ref{CKPAE}). \qed\\
So it is worth noting that $\tilde{L}_\phi$ is absolutely continuous, if $\frac{\kappa}{1+x^2}$ is integrable and  if $E$ is proper or  $1/\kappa$ is not proper.\\

 Let $\kappa$ be bounded. Then $L_\phi$, $W_\kappa$ are bounded,  $\tilde{L}_\phi=L_\phi$, $\tilde{W}_\kappa=W_\kappa$,  and 
  \begin{equation}\label{RBC} 
L_\phi=\,T^*\;W_\kappa\,T, \quad W_\kappa=T\,L_\phi \,T^*
\end{equation}
  Moreover, if  $E$ is proper or $1/\kappa$ is not proper, then even $L_\phi'=L_\phi$ by (\ref{CCKPAE}), whence
$$L_\phi=\,T^{-1}W_\kappa\,T$$
An example for the latter case   is $L_{\e^{-|x|}}$, which is injective. Generally,  if $\kappa > 0$  does not decrease too rapidly (so that $1/\kappa$ is proper) the kernel  of  $L_\kappa$  is not trivial. For an instructive  example see sec.\,\ref{LWHO}. The trivial example here is $\kappa=1_\R$. Note that $W(1_\R)=I_{L^2(\R_+)}$ and $L_{1_\R}$ is the orthogonal projection on $H_+$, and $T':H_+\to L^2(\R_+)$, $T'h:=P_+\mathcal{F}h$.
Examples for proper $E$ 
are the isomorphic pairs $W(1_E)\,\simeq\, \frac{1}{2}(I+H_E)$ with $H_E$ in  (\ref{THTE}),
which we like to write as 
\begin{equation}\label{IPWHSIO}
H_E\,=\, T^{-1}\,W_{21_E-1}\,T
\end{equation}
The case of the finite Hilbert transformation $H_{[-1,1]}$ is studied in detail in  \cite[sec.\,3.3]{C19}, \cite[(3.20)]{C19}.\\

In conclusion we remark  on the spectral representations of $\tilde{W}_\kappa$  and $\tilde{L}_\phi$  in (\ref{CKPAE}). By the spectral theorem in the multiplication operator version,
$\tilde{W}_\kappa$ is Hilbert space isomorphic to the multiplication operator $M(\varphi)$ on $L_\mu^2(\R)$ for some    Borel-measurable positive $\varphi:\R\to \R$ and finite Borel measure $\mu$. Let $V:L_\mu^2(\R)\to L^2(\R_+)$ be an isomorphism satisfying
\begin{equation}\label{SRWK}
\tilde{W}_\kappa = V \,M(\varphi)\,V^{-1}
\end{equation}
The  spectral measure  for $\tilde{W}_\kappa$ is given by $E_{\tilde{W}_\kappa}(\Delta)=VM(1_{\varphi^{-1}(\Delta)})V^{-1}$ for measurable $\Delta\subset \R$. 

\begin{Cor}\label{SPWF} 
Suppose the  representation \emph{(\ref{SRWK})} of $\tilde{W}_\kappa$.
 Then $W_0:=A^*VM(\frac{1}{\sqrt{\varphi}})$ with $A^*$ in \emph{(\ref{MLEM})} is closable, its closure $W$ is the Hilbert space isometry $T^*V$ with $\ran W=\big(\ker\tilde{L}_\phi\big)^\perp$, by which  $\tilde{L}_\phi=WM(\varphi)W^*$.
\end{Cor}\\
{\it Proof.} By (\ref{SRWK}) and (\ref{MLEM})(a), $W_0=A^*VM(\frac{1}{\sqrt{\varphi}})V^{-1}V=A^*\tilde{W_\kappa}^{-1/2}V=A^*|A^*|^{-1}V=T^*|_{\dom(|A^*|^{-1})}V$. For the last equality recall $A^*=T^*|A^*|$. Further note that $\dom(|A^*|^{-1})=\ran(|A^*|)$ is dense as $|A^*|=\tilde{W_\kappa}^{1/2}$ is self-adjoint injective. Since $T^*$ is bounded,
$W_0$ is closable and its closure $W$ equals $T^*V$. The remainder is obvious.\qed

\subsection{Example: Lalescu's operator}\label{LWHO} 
  Supposing  $\kappa>0$, $\kappa\in L^2$, and kernel $k\in L^1$, in \cite{P66} a spectral theory of $W_\kappa$ is proposed by  a reduction to a previously developed theory  for singular integral operators. Its application in \cite{P66} to Lalescu's  operator $W_\lambda$ with symbol 
\begin{equation}  
 \lambda(x):=\textrm{\small{$\frac{2}{1+x^2}$}}
  \end{equation}
   however does not produce the right normalization (\ref{NSFLO}) of the generalized eigenfunctions (\ref{SPLONN}). In establishing (\ref{NSFLO}) we get also the  result  (\ref{SOP}) on  orthogonal polynomials.\\
\hspace*{6mm}  
The diagonalization of Lalescu's operator $W_\lambda$ is achieved in sec.\,\ref{SRWL} and that of the associated singular integral operator $L_\lambda$ is derived  in sec.\,\ref{SRLL}. The generalized eigenfunctions  (\ref{GEFSIALO}) of the latter are no longer regular distributions.

\subsubsection{Spectral representation of $W_\lambda$}\label{SRWL}
Clearly,  the spectrum of $W_\lambda$ 
 lies in $[0,2]$. The integral kernel for $W_\lambda$ 
is $\operatorname{e}^{-|x|}$. 
 The obvious ansatz $a\e^{\alpha x}+b\e^{\beta x}$ for $u(x)$, $x>0$ in $s\,u(x)-\int_0^\infty \e^{-|x-y|} u(y)\d y=0$
yields for every  $s\in]0,2[$ the generalized  eigenfunction  (cf.\,\cite{WHP31})
\begin{equation}\label{SPLONN}
 q_s(x)=n(s)\left((\tau-\i )\e^{\i\tau x}+(\tau+\i)\e^{-\i\tau x}\right),\quad\tau:=\big(\textrm{\SMALL{$\frac{2}{s}$}}-1\big)^{1/2}
\end{equation}
 The claim is that there is the unique positive normalization constant 
\begin{equation}\label{NSFLO}
n(s)=(4\pi s\tau)^{-1/2}
\end{equation}
 such that $q(s,x):=q_s(x)$ is  the kernel for a Hilbert space isomorphism
\begin{equation}\label{HSIDLO} 
V:L^2(0,2)\to L^2(\R_+), \quad Vh=\operatorname{l.i.m.}\int_{\downarrow 0}^{\uparrow 2}q(s,\cdot)h(s)\d s
\end{equation}

Indeed, the assertion follows immediately from the following Lemma (\ref{HSIOP}) as $V=U^*\Gamma$ holds for the Hilbert space isomorphism 
\begin{equation}\label{IRL2S}
\Gamma:L^2(0,2)\to L^2(\R_+), \quad \Gamma h\,(y):=\textrm{\SMALL{$\frac{2\sqrt{y}}{1+y^2}$}} h\Big(\textrm{\SMALL{$\frac{2}{1+y^2}$}}\Big)
\end{equation}
 due to the change of variable $\gamma(y):=\frac{2}{1+y^2}$. 

\begin{Lem}\label{HSIOP}  $Ug(x):=\operatorname{l.i.m.}\int_{\downarrow 0}^{\uparrow \infty}u(x,y)g(y)\d y$ with  $$u(x,y)=(2\pi)^{-1/2}(1+x^2)^{-1/2}\big((x-\i )\e^{\i xy}+(x+\i )\e^{-\i xy}\big)$$ is unitary on $L^2(\R_+)$.
\end{Lem}\\
{\it Proof.}  Using the generating function $l_q(x)= (1-q)^{-1}\exp\big(-\frac{1}{2}\frac{1+q}{1-q}\, x\big)$, $|q|<1$ of the Laguerre functions  $l_n$, $n\ge 0$ one finds $\lambda_n(x):=Ul_n(x)=\sqrt{\frac{2}{\pi}}(1+x^2)^{-1/2}\Big(\frac{x+\i}{1+2\i x}\gamma(x)^n+\frac{x-\i}{1-2\i x}\gamma(x)^{-n}\Big)$ with $\gamma(x):=\frac{-1+2\i x}{1+2\i x}$. Note that $\int_{-\infty}^\infty(1+x^2)^{-1}\big(\frac{x+\operatorname{i}}{1+2\operatorname{i}x}\big)^2\gamma(x)^k\operatorname{d}x=0$ for $k\in\mathbb{Z}$, since the integral is a parametrization of 
the complex integration of $\frac{\operatorname{i}}{4}\frac{3z-\operatorname{i}}{3-z}z^k$ along the unit circle. Using  this it follows easily that $\lambda_n$, $n=0,1,2,\dots$ are orthonormal in $L^2(\R_+)$. It remains to show the completeness of $(\lambda_n)_n$.\\
\hspace*{6mm} Let $\varphi\in]0,\pi[$ satisfy $\operatorname{e}^{\operatorname{i}\varphi}=\gamma(x)$,  $x>0$. Put $\xi:=\cos\varphi$. Then $\lambda_n(x)=\sqrt{\frac{18(1-\xi)}{\pi(5-3\xi)}}\sqrt{1-\xi^2}\big(\frac{\sin(n+1)\varphi}{\sin \varphi}-\frac{1}{3}\frac{\sin n\varphi}{\sin \varphi}\big)=\sqrt{\frac{18
(1-\xi)}{\pi(5-3\xi)}}\sqrt{1-\xi^2}\big(U_n(\xi)-\frac{1}{3}U_{n-1}(\xi)\big)$, where $U_n$ denotes the  Chebyshev polynomial of the second kind of degree $n$ for $n=0,1,2,\dots$ and $U_{-1}:=0$ \cite[10.11(2)]{E53}. Therefore the Hilbert space isomorphism $\Gamma: L^2(\R_+)\to L^2(-1,1)$ due to the change of variable $\gamma(\xi):=\frac{1}{2}\sqrt{\frac{1+\xi}{1-\xi}}$, $\gamma'(\xi)=\frac{1}{2}\sqrt{\frac{1-\xi}{1+\xi}}(1-\xi)^{-2}$ maps $\lambda_n$ onto $\Gamma\lambda_n(\xi)=\sqrt{\frac{9}{\pi(5-3\xi)}}(1-\xi^2)^{1/4}Q_n(\xi)$, where $Q_n:=U_n-\frac{1}{3}U_{n-1}$ is a polynomial of degree $n$ for $n=0,1,2,\dots$  Thus (\ref{SOP}) is shown, which implies the assertion.\qed

\begin{Cor}\label{SOP} Let $Q_n:=U_n-\frac{1}{3}U_{n-1}$ , $n=0,1,2,\dots$ with $U_n$  the  Chebyshev polynomial of the second kind of degree $n$ and $U_{-1}:=0$. Then $(Q_n)_{n\ge0}$ is a sequence of orthonormal polynomials on ¢«$[-1,1]$ with respect to the weight function $\xi\mapsto\frac{9}{\pi}\frac{\sqrt{1-\xi^2}}{(5-3\xi)}$. It
obeys the recurrence  $Q_{n+1}=2\xi Q_n-Q_{n-1}$ for $Q_0=1$, $Q_1(\xi)=2\xi-\frac{1}{3}$.
\end{Cor}

\begin{The}\label{DLO} The spectral representation $W_\lambda=VM(\operatorname{id}_{[0,2]})V^{-1}$ holds. 
\end{The}\\
{\it Proof.} First note $|q(s,x)|\le \big(s(2-s)\big)^{-1/2}$. Let $h\in L^2(0,2)$ with support in $[\delta,2-\delta]$ for some $0<\delta<1$. Then $\big(VM(\operatorname{id}_{[0,2]})h\big)(x)=\int_\delta^{2-\delta}q(s,x)s\,h(s)\d s$. Substitute  $s\,q(s,x)=\int_0^\infty\e^{-|x-y|}q(s,y)\d y$. Obviously the integrations can be interchanged by Fubini's theorem yielding $\int_0^\infty\e^{-|x-y|}\Big(\int_\delta^{2-\delta}q(s,y)h(s)\d s\Big)\d y=\int_0^\infty\e^{-|x-y|}\,Vh(y)\d y=W_\lambda Vh(y)$. This implies the result.\qed

\subsubsection{Spectral representation of $L_\lambda$}\label{SRLL}
Let us first establish the kernel of $L_\lambda$. By (\ref{MLEM})(e) the latter equals that of $A=P_+\F M(\sqrt{\lambda})$, whence  $\ker(L_\lambda)=\{f\in L^2(\R):\sqrt{\lambda}f\in H_-\}$ by (\ref{CPEFIPS}). 
Put $u(x):=\frac{x+\i}{\sqrt{1+x^2}}$, $x\in \R$. Note that the multiplication operator $M(u)$ is unitary.  It follows
\begin{equation} 
\ker(L_\lambda)=M(u)H_-, \quad \ker(L_\lambda)^\perp=M(u)H_+
\end{equation}
Indeed, as to the less trivial implication let $f\in\ker(L_\lambda)$. Then $h:=\sqrt{\lambda}f\in H_-$, whence $(x-\i)h\in L^2(\R)$. Therefore $k:=2^{-1/2}(x-\i)h\in H_-$ using (\ref{AL1})(c), and $f=uk$ follows. \\
\hspace*{6mm} Hence the reduction $L'_\lambda$ of $L_\lambda$ on the orthogonal complement of its kernel is  self-adjoint bounded on $M(u)H_+$ with spectrum $[0,2]$. The following computations are valid in a  distributional sense.  So generalized eigenfunctions $Q_s$ for $L_\lambda$  are given by $A^*q_s = M(\sqrt{\lambda})\F^{-1}P_+^*q_s$ since $L_\lambda A^*q_s=A^*A\,A^*q_s=sA^*q_s$. Recall $\int_{-\infty}^\infty 1_{]0,\infty[} \e^{\i xy}\d y=\frac{\i}{x}+\pi\delta(x)$, see e.g. \cite[Table of Fourier Transforms 1.23]{GS64}. Hence one yields the not regular distributions on $\R$
\begin{equation}\label{GEFSIALO}
Q_s(x) =\textrm{\SMALL{$\frac{n(s)}{\sqrt{s}}$}}\textrm{\SMALL{$\sqrt{\frac{2}{1+x^2}}\;\frac{1}{\sqrt{2\pi}}$}}\left((\tau-\i)\Big(\frac{\i}{x+\tau}+\pi\delta(x+\tau)\Big)+(\tau+\i)\Big(\frac{\i}{x-\tau}+\pi\delta(x-\tau)\Big)\right)
\end{equation}
The claim is that $\frac{n(s)}{\sqrt{s}}$ is the unique positive normalization constant such that $Q(s,x):=Q_s(x)$ is  the kernel for an isometry  $W$ on $L^2(0,2)$ in $L^2(\R)$ satisfying
\begin{equation}\label{HSIDLO} 
 Wh=\int_ 0^2 Q(s,\cdot)h(s)\d s
\end{equation} 
for test functions $h$. The additional factor $\frac{1}{\sqrt{s}}$ regarding the normalization constant of $Q_s$ corresponds to the factor $\frac{1}{\sqrt{\varphi}}$ for $W_0$ in (\ref{SPWF}) and  is suggested heuristically by $\langle A^*q_s, A^*q_s\rangle=\langle q_s,AA^*q_s\rangle=s\langle q_s,q_s\rangle$. 
For the proof recall the isomorphism 
$\Gamma$ (\ref{IRL2S}) and  the Hilbert transformation $H$ (\ref{SRHT}). For $g:\R_+\to \C$ let $g_{oe}$ denote the odd extension $g_{oe}(x)=-g(-x)$ for $x\le 0$ of $g$.

\begin{The}  The integral operator \emph{(\ref{HSIDLO})} determines a Hilbert space isometry $W: L^2(0,2)\to L^2(\R)$ with $\ran W=M(u)H_+$. It satisfies 
 $Wh=M(u)\,\frac{1}{2}\big(I+H\big)(\Gamma h)_{oe}$ for $h\in L^2(0,2)$ and yields the representation 
 $$L_\lambda=WM(\operatorname{id}_{[0,2]})W^*$$
\end{The}\\
{\it Proof.} By the change of variable $t=2/(1+s^2)$ for the integration  in (\ref{HSIDLO}) one obtains $\sqrt{1+x^2}\,Wh(x)=\frac{1}{2\pi}\int_0^\infty \left((t-\i)\Big(\frac{\i}{x+t}+\pi\delta(x+t)\Big)+(t+\i)\Big(\frac{\i}{x-t}+\pi\delta(x-t)\Big)\right)\Gamma h(t)\d t$. The integral is easily done yielding ($\star$)  $Wh=M(u)\,\frac{1}{2}\big(I+H\big)(\Gamma h)_{oe}$  for test functions $h$. Check $\langle Hf,f\rangle=0$ if $f\in L^2(\R)$ is odd. Therefore $\norm{Wh}^2=\norm{\frac{1}{2}\big(I+H\big)(\Gamma h)_{oe}}^2=\langle\frac{1}{2}\big(I+H\big)(\Gamma h)_{oe},(\Gamma h)_{oe}\rangle=\frac{1}{2}\norm{(\Gamma h)_{oe}}^2=\norm{h}^2$. So $W$ is an isometry, and ($\star$) extends to all $h\in L^2(0,2)$. Obviously $Wh\in M(u)H_+$. Moreover, if $k\in H_+$ then $k=\frac{1}{2}\big(I+H\big)(k-\check{k})$ as $\check{k}\in H_-$. Since $k-\check{k}$ is odd, this implies $\ran W=M(u)H_+$.\\
\hspace*{6mm}
The  spectral representation of $L_\lambda$ follows from $WM(\operatorname{id}_{[0,2]})h=L_\lambda Wh$ for test functions $h$. The latter is shown along the lines of the proof of (\ref{DLO}). We omit the technicalities.\qed

\end{document}